\documentclass[12pt,a4paper]{article}
\usepackage[T2A]{fontenc}
\usepackage[cp1251]{inputenc}
\usepackage{epsfig,amssymb,epic,eepic,color, graphicx}
\usepackage{fancyhdr}
\usepackage{amsthm}
\usepackage{amsmath}
\usepackage{indentfirst}
\usepackage{float}
\usepackage{setspace}
\usepackage{float}
\usepackage{fancyhdr}
\newtheorem{theorem}{Theorem}

\newtheorem{collary}{Collary}[section]

\newtheorem{lemma}{Lemma}[section]
\newtheorem{proposition}{Proposition}[section]
\newenvironment{demo}{{\bf Proof: }}{\hfill $\square$ \medskip}
\begin{document}

UDC: 517.938.5


\begin{center}{\bf Elena V. Nozdrinova and Olga V. Pochinka}
\\[1 true cm]
{\large\bf Stable arcs connecting polar cascades on a torus}
\\[1 true cm]\end{center}
E.V. Nozdrinova: maati@mail.ru,    HSE, B. Pecherskaya
street, 25/12, Nizhny Novgorod 603150, Russia; ORCID: 0000-0001-5209-377X
\\[1 true cm]
O.V. Pochinka: olga-pochinka@yandex.ru,   HSE, B. Pecherskaya street, 25/12, Nizhny Novgorod 603150, Russia; ORCID: 0000-0002-6587-5305
\\[1 true cm]
{\bf Abstract}
In this paper, we obtain a solution to the 33rd Palis-Pugh problem for polar gradient-like diffeomorphisms on a two-dimensional torus, under the assumption that all non-wandering points are fixed and have a positive orientation type.


{\it Acknowledgements.}  
Study of the dynamics of diffeomorphisms of the class under consideration is supported by RSF (Grant No. 17-11-01041), the construction of a stable arc is supported by Laboratory of Dynamical Systems and Applications NRU HSE, of the Ministry of science and higher education of the RF grant ag. No 075-15-2019-1931 

\section{Introduction and formulation of results} 
The problem of the existence of an arc with no more than a countable (finite) number of bifurcations connecting structurally stable systems (Morse-Smale systems) on manifolds is on the list of fifty Palis-Pugh problems \cite{PaPu} under number 33.

In 1976, S. Newhouse, J. Palis, F. Takens \cite{NPT1976} introduced the concept of a stable arc connecting two structurally stable systems on a manifold. Such an arc does not change its quality properties with a small perturbation. In the same year, S. Newhouse and M. Peixoto \cite{NP} proved the existence of a simple arc (containing only elementary bifurcations) between any two Morse-Smale flows. It follows from the result of G. Fleitas \cite{Fl} that a simple arc constructed by Newhouse and Peixoto can always be replaced by a stable one \cite{NPT}. For Morse-Smale diffeomorphisms given on manifolds of any dimension, examples of systems that cannot be connected by a stable arc are known. 

Obstruction appear already for orientation-preserving diffeomorphisms of the circle $S^1$, which are connected by a stable arc only if the rotation numbers coincide \cite {No}.

Beginning with dimension two, additional obstruction appear to the existence of stable arcs between isotopic diffeomorphisms. They are associated with the existence of periodic points \cite {Bla}, \cite {NoPo20}, heteroclinic intersections \cite {Ma}, wild embeddings of separatrices \cite {GP}, etc.

On the 6-dimensional sphere, examples of source-sink diffeomorphisms are known that are not connected by any smooth arc \cite {BoGrPo}, which, in fact, became the source for constructing different smooth structures on a sphere of dimension 7. For $n=2,3$ the non-trivial fact of the existence of an arc without bifurcations between two source-sink diffeomorphisms was established in \cite{BoGrPo}, \cite{NoPo}.   

A natural generalization of source-sink systems are {\it polar diffeomorphisms} -- gradient-like diffeomorphisms with a unique source and a unique sink. It follows from Morse theory that such diffeomorphisms exist on any manifolds. 

In this paper, we consider the class $G$ of polar gradient-like diffeomorphisms on the two-dimensional torus $\mathbb T^2$, under the assumption that all non-wandering points are fixed and have positive orientation type. In chapter \ref{dp} it is established that any diffeomorphism $ f \in G $ has exactly two saddle points and is isotopic to the identity. Moreover, all diffeomorphisms of the class under consideration are pairwise topologically conjugate (see, for example, \cite {BeGr}, \cite {GrKaPo}). Moreover, the closures of stable (unstable) manifolds of saddle points of different diffeomorphisms can belong to different homotopy classes of closed curves on the torus. Therefore that in the general case there is no arc without bifurcations connecting two diffeomorphisms of the class under consideration. 

The main result of this work is the proof of the following theorem.

\begin{theorem}\label{is} 
Any diffeomorphisms $f, f'\in G$ can be connected by a stable arc with a finite number of saddle-node bifurcations.
\end{theorem}

\section{Diffeomorphisms of class $G$}\label{dp}

\subsection{On gradient-like flows on surfaces}
In this section, we establish the basic dynamical properties of diffeomorphisms $f:\mathbb T^2 \to \mathbb T^2$ from the class $G$.

Recall that a diffeomorphism $f$ is {\it gradient-like} if its non-wandering set $ \Omega_f $ consists of a finite number of hyperbolic points and the invariant manifolds of different saddle points do not intersect. A diffeomorphism $f$ is called {\it polar} if the set $ \Omega_f $ contains exactly two nodal points, namely, one sink and one source.

Fix system of generators of fundamental group of torus $\mathbb T^2=\mathbb S^1\times\mathbb S^1$: $$a=\mathbb S^1\times\{0\}=<1,0>;\,\,b=\{0\}\times\mathbb S^1=<0,1>.$$ 
Recall that the {\it algebraic torus automorphism} $\widehat L: \mathbb T^2\to \mathbb T^2$, $\mathbb T^2=\mathbb R^2/\mathbb Z^2$ is called the diffeomorphism defined by the matrix  $\,{L}=
\left(\begin{array}{cc}
 \alpha & \beta \\
 \gamma & \delta
\end{array} \right)$, belonging to the set  $GL(2,\mathbb Z)$ {\it unimodular matrices} -- integer matrices with determinant $\pm 1$. That is $$\widehat L(x,y)=(\alpha x+\beta y,\gamma x+\delta y) \pmod 1.$$

\begin{theorem}\label{dynG} Any diffeomorphism $f \in G$ has the following properties:
\begin{enumerate}
\item Non-wandering set $\Omega_f$ of the diffeomorphism $ f $ consists of exactly four fixed hyperbolic points: the sink $ \omega_f $, the source $ \alpha_f $, and the saddles $\sigma_f^1, \sigma_f^2$, the closures of invariant manifolds of which are closed curves: $$c^{s1}_f=cl\, W^s_{\sigma_f^1}=W^s_{\sigma_f^1}\cup\alpha_f,\,\, c^{u1}_f=cl\, W^u_{\sigma_f^1}=W^u_{\sigma_f^1}\cup\omega_f,$$ $$c^{s2}_f=cl\, W^s_{\sigma_f^2}=W^s_{\sigma_f^2}\cup\alpha_f,\,\, c^{u2}_f=cl\, W^u_{\sigma_f^2}=W^u_{\sigma_f^2}\cup\omega_f.$$
\item There is only one choice of saddle points numbering $\sigma_f^1,\,\sigma_f^2$ and the orientation of the closures of their invariant manifolds such that the curves $c^{s1}_f,\,c^{u2}_f$ are of homotopy type $<\mu_f^1,\nu_f^1>$ and the curves  $c^{s2}_f,\,c^{u1}_f$ are of homotopy type $<\mu_f^2,\nu_f^2>$ in the basis $a,b$, also 
$J_f=\begin{pmatrix}
  \mu^1_f& \mu_f^2\\
  \nu_f^1& \nu_f^2
\end{pmatrix}$ is unimodular matrix with the following properties:

a) $\mu^1_f\geq\mu^2_f\geq 0$,

b) $\nu^1_f>\nu^2_f$, if $\mu^1_f=\mu^2_f$,

c) $\nu^2_f=1$, if $\mu_f^2=0$.

\item The diffeomorphism $f$ is isotopic to the identity map.
\end{enumerate}
\end{theorem}
\begin{demo} Let $f\in G$. We prove all the points of the theorem sequentially.

1. We denote by $ k_f^0, k_f^1, k_f^2 $ the number of sinks, saddles and sources of the diffeomorphism $f$ respectively.  According to \cite {Pi}, there is a Morse function on the torus $\mathbb T^2$, the set of critical points of which coincides with the set $\Omega_f$ and the indices of the critical points coincide with the dimensions of the unstable manifolds of non-wandering points of the diffeomorphism $ f $. Then it follows from Morse inequalities (see, for example, \cite {Mil}) that
$$k_f^0-k_f^1+k_f^2=0.$$
Since the diffeomorphism $f$  is polar, then $k_f^0=k_f^2=1$ and therefore $k_f^1=2$. Thus, the non-wandering set $\Omega_f$  of the diffeomorphism $f$ consists of exactly four fixed hyperbolic points: the sink $\omega_f$, the source $\alpha_f$ and the saddles $\sigma_f^1, \sigma_f^2$.

Since the invariant manifolds of different saddle points of the diffeomorphism $f$ do not intersect, then, according to \cite[proposition  2.1.3]{GrPo},  
$$cl(W^u_{\sigma_f^1})\setminus W^u_{\sigma_f^1}=cl(W^u_{\sigma_f^2})\setminus W^u_{\sigma_f^2}=\omega_f,$$
$$cl(W^s_{\sigma_f^1})\setminus W^s_{\sigma_f^1}=cl(W^s_{\sigma_f^2})\setminus W^s_{\sigma_f^2}=\alpha_f$$ and sets  $c^{u1}_f=cl\, W^u_{\sigma_f^1},\,\, c^{s1}_f=cl\, W^s_{\sigma_f^1},\,\,c^{u2}_f=cl\, W^u_{\sigma_f^2},\,\, c^{s2}_f=cl\, W^s_{\sigma_f^2}$ are homeomorphic to circles.

2. Denote by $<\mu_f^1,\nu_f^1>,\,\,<\mu_f^2,\nu_f^2>$ homotopy types of oriented curves
$c^{s1}_f,\,c^{u1}_f$, respectively, in the basis $a,b$. Let $J_f=\begin{pmatrix}
  \mu^1_f& \mu_f^2\\
  \nu_f^1& \nu_f^2
\end{pmatrix}.$ The previous item implies,  that the closed curves $c^{s1}_f,\, c^{u1}_2$ have a unique transverse intersection point $\sigma_f^1$. Then the modulus of the intersection index of these curves is one. Since this index coincides with the determinant of the matrix $J_f$ (see, for example, \cite[exercise 7, pp. 28]{Rolf}), the matrix $J_f$ is unimodular. Сonsequently, the curves $c^{s1}_f,\, c^{u1}_f$ are not homotopic to zero. Without loss of generality, suppose the orientation on the curves was chosen so that $\mu^i_f\geq 0$ and $\nu^i_f=1$, if $\mu^i_f= 0$ (these conditions uniquely orient the curves).

Since the diffeomorphism $f$ is gradient-like, then the curves $c^{u1}_f,\,c^{s2}_f$ and $c^{u2}_f,\,c^{s1}_f$ are pairwise disjoint. Then the oriented curves $c^{u2}_f,\,c^{s2}_f$ also have homotopy types $<\mu_f^1,\nu_f^1>$ and $<\mu_f^2,\nu_f^2>$, respectively, in the basis $a,b$ (see, for example, \cite[Theorem 13,  pp. 25]{Rolf}). Without loss of generality, suppose the numbering of saddle points is chosen so that  $\mu^1_f\geq\mu^2_f$ and, if  $\mu^1_f=\mu^2_f$, then $\nu^1_f>\nu^2_f$ (these conditions determine the numbering of saddle points uniquely).

3. The diffeomorphism $f$ induces an isomorphism $f_* : \pi_1(\mathbb T^2) \to \pi_1(\mathbb T^2)$ in the fundamental group $\pi_1(\mathbb T^2)$ of the torus, isomorphic to the abelian group $\mathbb Z^2$.  Then the isomorphism $f_*$ is determined by the unimodular integer matrix uniquely $L_f=\begin{pmatrix}
  \alpha & \beta\\
  \gamma & \delta
\end{pmatrix},$ taking basis $a, b$ into basis $<\alpha,\gamma>,\,<\beta,\delta>$.  Diffeomorphism $f$ isotopic to the identical iff $L_f=E$, where $E=\begin{pmatrix}
  1& 0\\
  0& 1
\end{pmatrix}$ (see, for example, \cite[lemma 3, pp. 26]{Rolf}). Let us show that $L_f=E$ for $f\in G$. 

Let $h=\widehat J^{-1}_f f \widehat J_f:\mathbb T^2\to\mathbb T^2$. By construction, the diffeomorphism $h$ is smoothly conjugate to the diffeomorphism $f$. Moreover, the curves $c^{s1}_h, c^{u2}_h$ and $c^{s2}_h, c^{u1}_h$ are of homotopy types $<1,0>$ and  $<0,1>$, respectively. Since all these curves are $h$-invariant, then $L_h=E$. Since $f=\widehat J_f h\widehat J_f^{-1}$, then $L_f=J_f L_hJ_f^{-1}=J_f EJ_f^{-1}=E$.
\end{demo}

\section{Construction of a stable arc between diffeomorphisms of the class $G$} 

In this section, we outline the proof of theorem \ref{is} with references to the statements that will be proved in the following sections. Let us first give the necessary definitions.

Consider a 1-parametric family of diffeomorphisms ({\it an arc}) $\varphi_t: M \to M, t \in [0,1]$. An arc $\varphi_t$ is called {\it smooth},  if map $F:M\times[0,1]\to M$,  defined by the formula $F(x,t)=\varphi_t(x)$ is smooth.

The smooth arc $\varphi_t$ is called a {\it smooth} product of the smooth arcs $\varphi^1_t$ and $\varphi^2_t$ such that $\varphi^1_1=\varphi^2_0$, if 
$\varphi_t =\begin{cases}\varphi^1 _ {2\tau (t)},~ 0\leqslant t\leqslant\frac {1} {2},\\
\varphi^2_{2\tau(t)-1}, \frac {1} {2} \leqslant t\leqslant 1,\end{cases}$ where $\tau: [0,1] \to [0,1] $ is a smooth monotone map such that $\tau (t) =0$
for $0\leqslant t\leqslant\frac {1} {3} $ and $\tau (t)
=1$ for $\frac {2} {3} \leqslant t\leqslant 1$. We will write $\varphi_t=\varphi^1_t*\varphi^2_t$. 

Following \cite{NPT}, an arc $\varphi_t$ is called {\it stable} if it is an inner point of the equivalence class with respect to the following relation: two arcs $\varphi_t$, $\varphi'_t$ are called {\it conjugate} if there
are homeomorphisms $ h: [0,1] \to [0,1 ], \, H_t: M \to M $ such that  $ H_t \varphi_t = \varphi' _ { {h (t )}} H_t, t \in [0, 1] $ and $ H_t $ continuously depend on $ t $.

In \cite{NPT} also established that the arc $\{\varphi_t\}$, consisting of diffeomorphisms with a finite limit set, is stable iff all its points are structurally stable diffeomorphisms with the exception of a finite number of bifurcation points, $\varphi_{b_i},i=1,\dots,q$ such that $\varphi_{b_i}$: 

1) has no cycles;

2) has a unique non-hyperbolic periodic orbit, which is a non-critical saddle-node or flip;

3) the invariant manifolds of all periodic points of the diffeomorphism$\varphi_{b_i}$  intersect transversally;

4) the transition through $\varphi_{b_i}$ is a generically unfolded saddle-node or period doubling bifurcation, wherein the saddle-node point is non-critical.

Recall the definition of unfolding generically arc $\varphi_t$ through the saddle-node or flip. We give the definition for a fixed non-hyperbolic point, in the case  when it has a period $k>1$, a similar definition is given for the arc $\varphi^k_t$.   

\begin{figure}
\centerline{\includegraphics[width=14 cm]{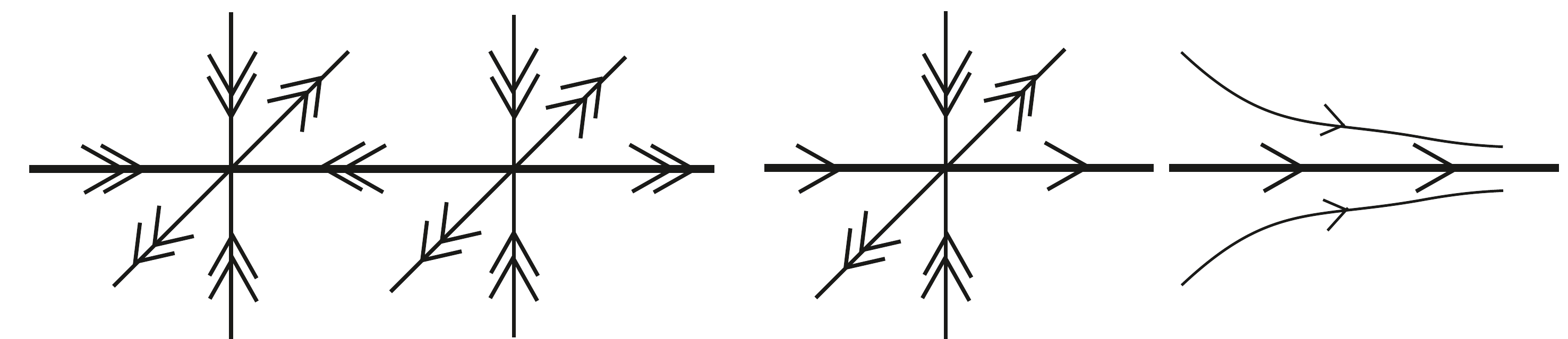}}
\caption{\small Saddle-node bifurcation}\label{saddle-node}
\end{figure}

An arc $\{\varphi_t\}\in \mathcal{Q}$ {\it unfolds generically} through a saddle-node bifurcation $\varphi_{b_i}$ (Fig. \ref{saddle-node}), if in some neighborhood of the nonhyperbolic point $(p,b_i)$ the arc $\varphi_t$ is conjugate to
$$\tilde\varphi_{\tilde t}(x_1,x_2,\dots,x_{1+n_u},x_{2+n_u},\dots,x_{n})=$$ $$\left(x_1+\frac {x_1^2}{2}+\tilde t,\pm 2x_2,\dots,\pm 2 x_{1+n_u},\frac{\pm x_{2+n_u}}{2},\dots,\frac{\pm x_{n}}{2}\right),$$  where $(x_1,\dots,x_n)\in\mathbb R^n,\,|x_i|< 1/2,\,|\tilde t|<1/10$.  

In the local coordinates $(x_1,\dots,x_n,\tilde t)$ the bifurcation occurs at time $\tilde t=0$ and the origin $O\in\mathbb R^n$ is a saddle-node point. The axis $Ox_1$ is called a {\it central manifold} $W^c_O$, the half-space   $\{(x_1,x_2, \dots, x_{n})\in\mathbb R^n:\,x_1\geq 0,\,x_{2+n_u}=\dots=x_n=0\}$  is the {\it unstable manifold} $W^u_O$, half-space $\{(x_1,x_2, \dots, x_{n})\in\mathbb R^n:\,x_1\leq 0,\,x_{2}=\dots=x_{1+n_u}=0\}$ is the {\it stable manifold} $W^s_O$  of the point $O$. 

\begin{figure}[h]
\centerline{\includegraphics[width=6cm]{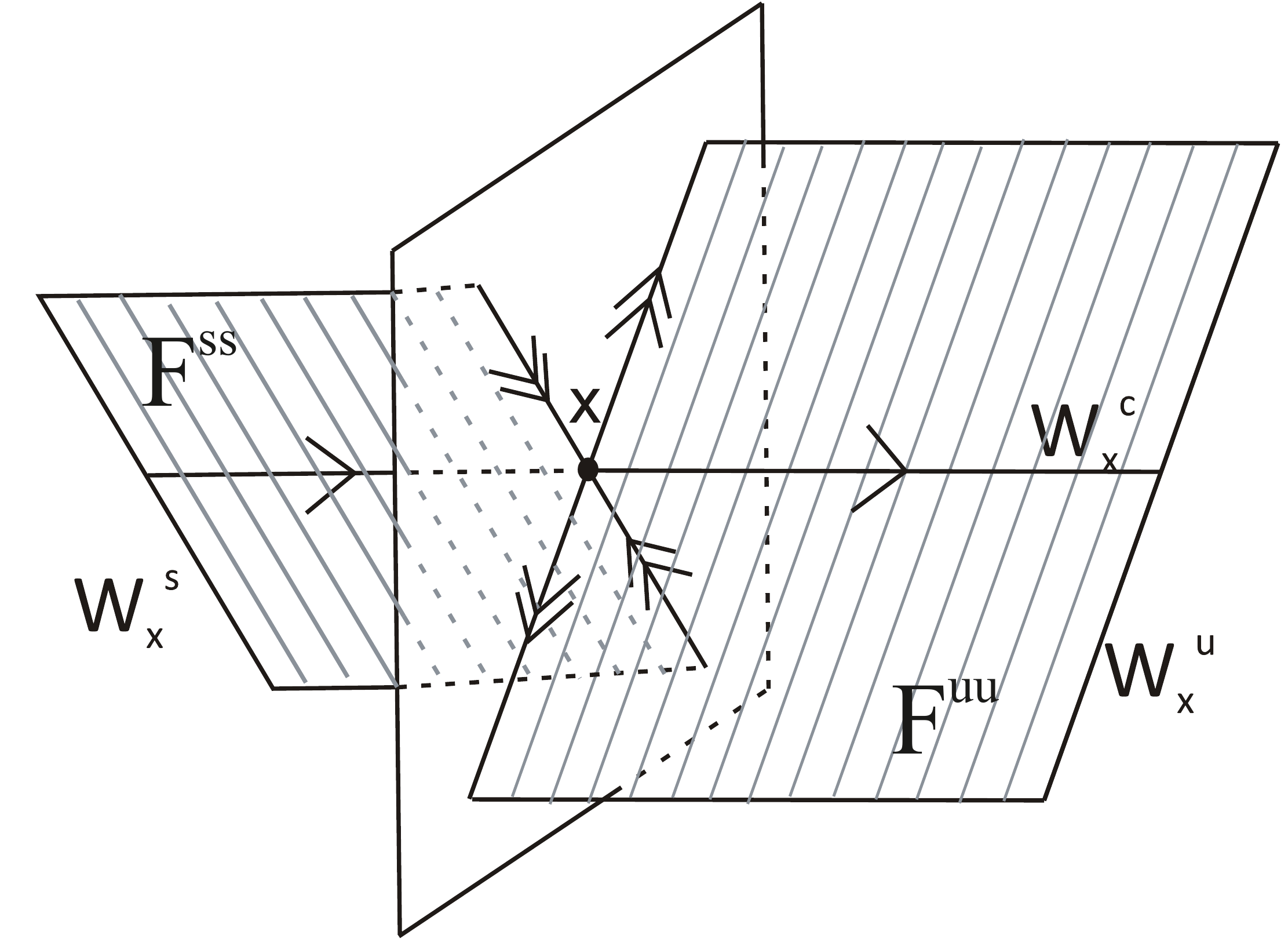}}
\caption{\small Strongly stable and unstable foliations}\label{NP}
\end{figure}

If $p$ is a saddle-nodal point of the diffeomorphism $\varphi_{b_i}$, then there exists a unique $\varphi_{b_i}$- invariant foliation $F^{ss}_p$ with smooth leaves such that $\partial W^s_p$ is a leave of this foliation \cite{HPS}. $F^{ss}_p$ is called a {\it strongly stable foliation} (Fig. \ref{NP}).
A similar {\it strongly unstable foliation} is denoted by $F^{uu}_p$. A point $p$ is called {\it $s$-critical}, if there exists some hyperbolic periodic point $q$ such that $W^u_q$ non-transversally intersect some leaf of the foliation $F^{ss}_p$; {\it $u$-criticality} is defined similarly. Point $p$ is called

- {\it semi-critical} if it is either  $s$- or $u$-critical;

- {\it bi-critical} if it is $s$- and $u$-critical;

- {\it non-critical} if it is not semi-critical\footnote{For the first time, the effect of arc instability in a neighborhood of a non-critical saddle- was discovered in 1974 by V. Afraimovich and L. Shilnikov \cite{ASh1}, \cite{ASh2}. The existence of invariant foliations $F ^{ss}_p, \, F^{uu}_p$ was also proved earlier in the works of V. Lukyanov and L. Shilnikov\cite{LSh}.}.

Let ${f},f'\in{G}$. Let us prove that the diffeomorphisms $f,f'$ are connected by a stable arc $\varphi_{t}: \mathbb T^{2} \rightarrow \mathbb T^{2}, t\in [0,1]$, whose diffeomorphisms are gradient-like except for a finite number of generically unfolding non-critical saddle-node bifurcations.

\begin{demo} In section \ref{mod} for any unimodular matrix $J=\begin{pmatrix}
  \mu^1& \mu^2\\
  \nu^1& \nu^2
\end{pmatrix}$ such that $\mu^1\geq\mu^2\geq0$ and $\nu^1>\nu^2$, if $\mu^1=\mu^2$, we construct the a model diffeomorphism ${f}_J\in{G}$, for which $J_{f_J}=J$. According to lemma \ref{fA}, every diffeomorphism $f\in G$  can be connected with the model diffeomorphism $f_{J_f}$ by an arc without bifurcations $H_{f,t}$. According to lemma \ref{AA}, the diffeomorphism $f_J$ can be joined by an arc $H_{J,t}$ with a finite number of generically unfolding non-critical saddle-node bifurcations with diffeomorphism  $f_{0}$. Then the desired arc $\varphi_t$ has the form $$\varphi_t=H_{f,t}*H_{J_f,t}*H_{J_{f'},1-t}*H_{f',1-t}.$$  
\end{demo}

\section{Construction of model diffeomorphisms in the class $G$}\label{mod}
In this section for any unimodular matrix $J=\begin{pmatrix}
  \mu^1& \mu^2\\
  \nu^1& \nu^2
\end{pmatrix}$ such that $\mu^1\geq\mu^2\geq0$ and $\nu^1>\nu^2$ if $\mu^1=\mu^2$ we construct a model diffeomorphism ${f}_J\in{G}$, for which $J_{f_J}=J$.

The simplest example of a diffeomorphism from the class $G$ is the direct product of two copies of a source-sink diffeomorphism on the circle $\mathbb S^1$, denote it by $f_0$. First, we construct a source-sink diffeomorphism on the circle. In order to do this consider the map $\bar{F}_0:\mathbb R\to\mathbb R$, given by the formula:
\begin{center}
$\bar{F}_0(x)=x-\frac{1}{4\pi} sin \left(2\pi\left (x-\frac14\right)\right).$
\end{center}

\begin{figure}[h!]
\centerline{\includegraphics[width=10 cm]{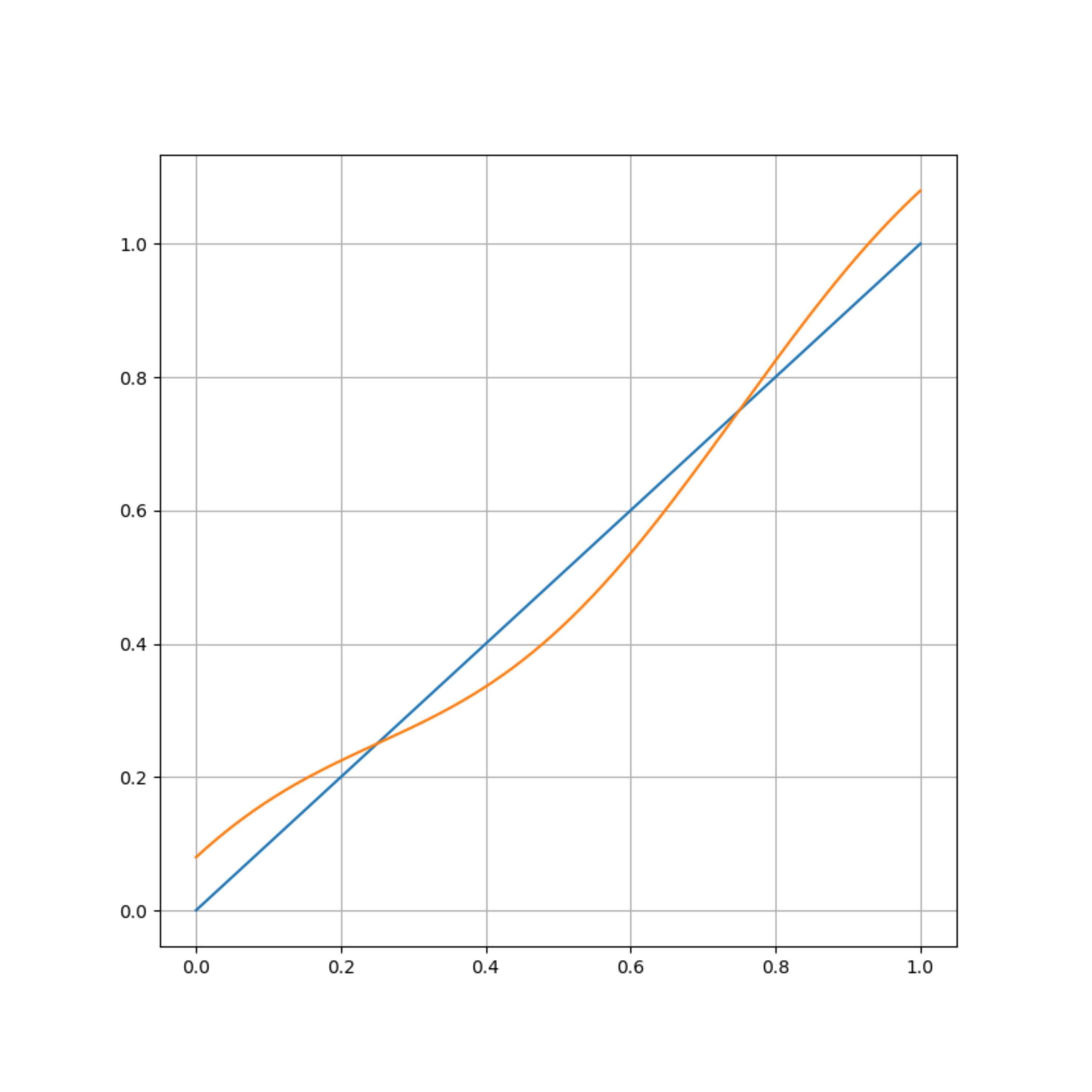}}
\caption{\small Graph of the map $\bar F_0$}\label{phi0}
\end{figure}

By construction $x=\frac14$ and $x=\frac34$ -- are fixed points of the map $\bar{F}_0$ on the segment $[0,1]$ (Fig. \ref{phi0}). 

Consider the projection $\pi: \mathbb{R}\to \mathbb{S}^1$ given by the formula $\pi(x)=e^{2\pi i x}$. As $\bar{F}_0$ is strictly increasing and satisfy the condition $\bar{F}_0(x+1)=\bar{F}_0(x)+1$, there is a diffeomorphism projecting it to the circle $$F_0=\pi\bar  F_0 \pi^{-1}:\mathbb S^1\to\mathbb S^1.$$ 

By construction, the diffeomorphism $F_0$ has a fixed hyperbolic sink at the point $N=\pi\left(\frac14\right)$ and a fixed hyperbolic source at the point $S=\pi\left(\frac34\right)$. 

\begin{figure}[h]
\centerline{\includegraphics[width=13 cm]{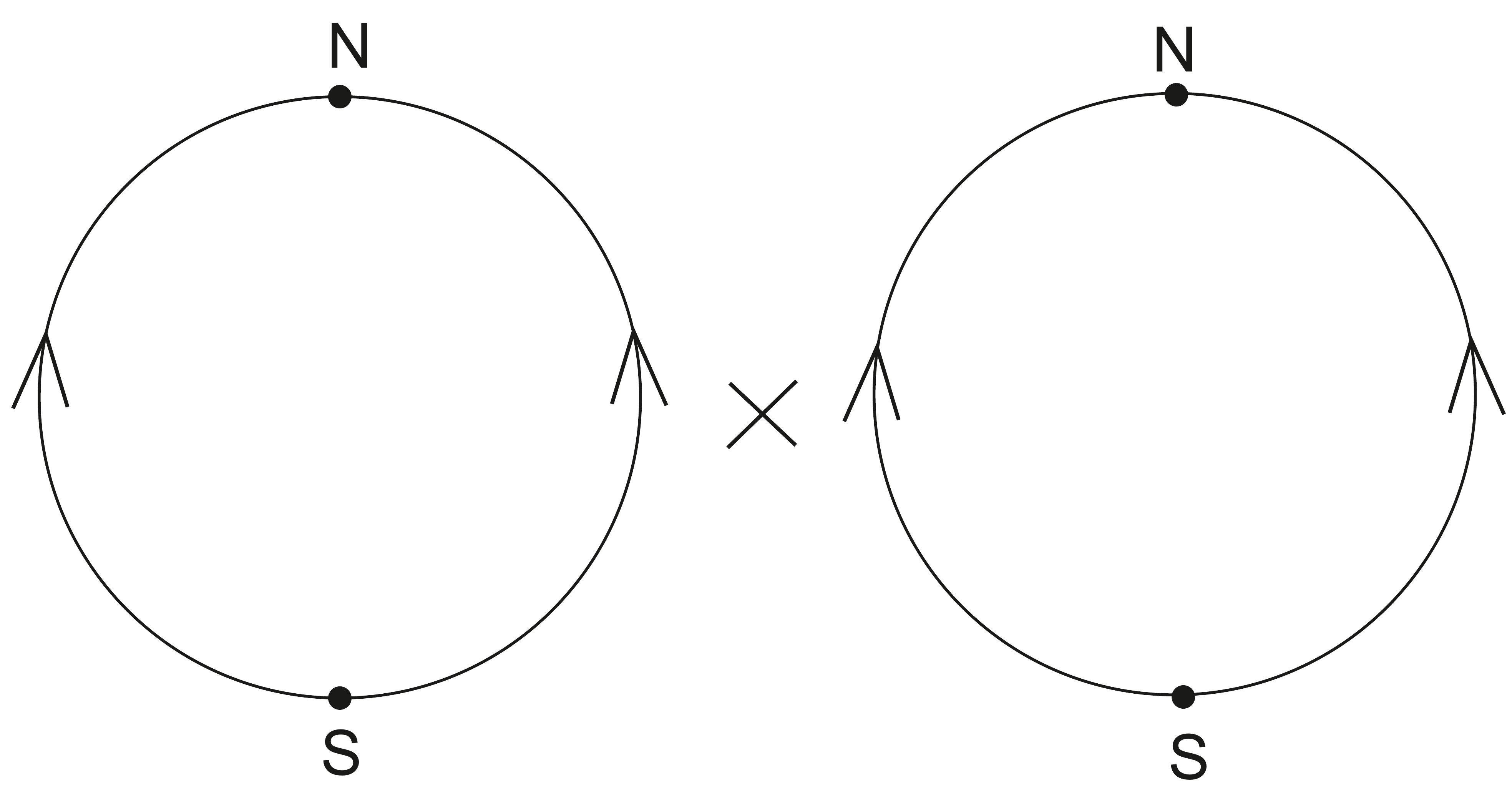}}
\caption{\small Cartesian square of the diffeomorphism $F_0$}\label{okr}
\end{figure}

Define the diffeomorphism $f_0:\mathbb T^2\to\mathbb T^2$ by the formula (Fig. \ref{okr}) $$f_0(z,w)=(F_0(z),F_0(w)),\,z,w\in\mathbb S^1.$$

\begin{figure}[h!]
\centerline{\includegraphics[width=5 cm]{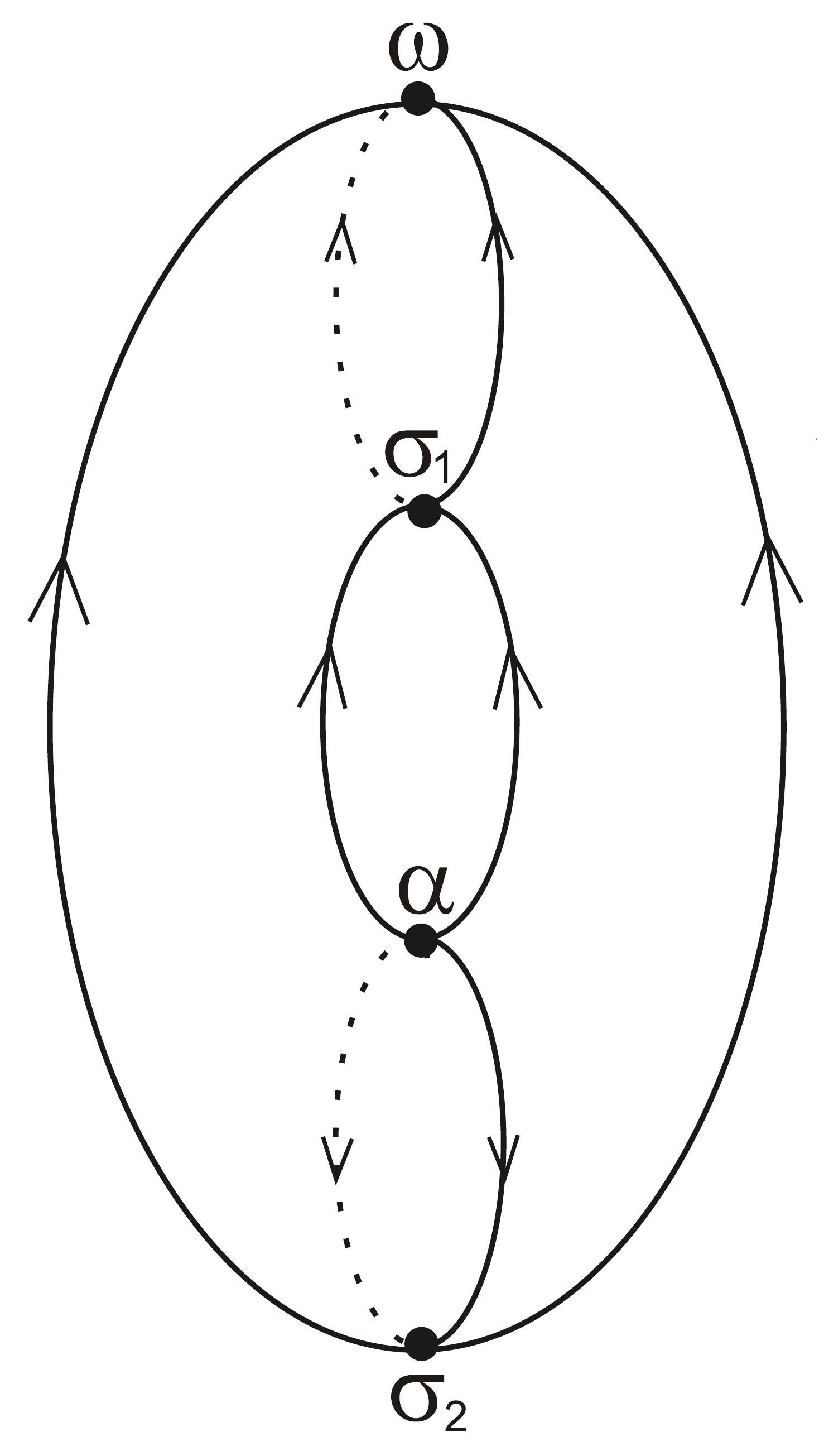}}
\caption{\small Diffeomorphism $f_0$.}\label{f0}
\end{figure}

By construction the diffeomorphism $f_0$ contain a fixed hyperbolic sink at the point $\omega=(N,N)$, hyperbolic source $\alpha=(S,S)$ and has two saddle points $\sigma_1=(N,S), \sigma_2=(S,N)$ (Fig. \ref{f0}). Moreover, the closures of their invariant manifolds lie in the classes of generators $a,b$. Exactly, $$c^{s1}_{f_0}=cl\,W^s_{\sigma_1}=\mathbb S^1\times\{S\},\,c^{u1}_{f_0}= cl\,W^u_{\sigma_1}=\{N\}\times\mathbb S^1,$$ $$c^{s2}_{f_0}= 
cl\,W^s_{\sigma_2}=\{S\}\times\mathbb S^1,\,c^{u2}_{f_0}=
cl\,W^u_{\sigma_2}=\mathbb S^1\times\{N\}.$$

Let $f_J=\widehat J f_0\widehat J^{-1}$. We will call the diffeomorphism $f_J$ {\it model diffeomorphism}.
By construction $f_E=f_0$.

\section{Construction the arc $H_{f,t}$}\label{A-A}
\subsection{Outline of the construction}
\begin{lemma}\label{fA} Every diffeomorphism $f\in G$ is connected by an arc without bifurcations $H_{f,t}$ with the diffeomorphism $f_{J_f}$.
\end{lemma}
\begin{demo} We give a scheme of the proof of Lemma \ref {fA} with links to statements that will be proved below.

Consider two possible cases: 1) $J_f=E$, 2) $J_f\neq E$.

1) According to the lemma \ref{PR4.3}, we assume that the diffeomorphism $f$  in a neighborhood of the sink $\omega_f$ has a local chart $(U_{\omega_f}, \psi_{\omega_f}), \psi_{\omega_f} : U_{\omega_f} \to \mathbb R^2$ such that the diffeomorphism $g=\psi_{\omega_f} f \psi_{\omega_f}^{-1}: \mathbb R^2 \to \mathbb R^2$ has the form $g(x,y)=\left(\frac{x}{2},\frac{y}{2}\right)$. Moreover, according to  the lemma \ref{L5}, we can assume that $$\psi_{\omega_f}(c^{u2}_{f}\cap U_{\omega_f})\subset Ox,\,\psi_{\omega_f}(c^{u1}_{f}\cap U_{\omega_f})\subset Oy,$$ that is, the curves $c^{u2}_{f},\,c^{u1}_{f}$ are smooth. Since they lie in the same homotopy class with curves $c^{u2}_{f_0}=\mathbb S^1\times\{N\},\,c^{u1}_{f_0}=\{N\}\times\mathbb S^1$, then, by Proposition \ref{Rol}, there exists a diffeomorphism $\xi:\mathbb T^2\to\mathbb T^2$ smoothly isotopic to the identity such that $$\xi(c^{u1}_{f})=c^{u1}_{f_0},\,\xi(c^{u2}_{f})=c^{u2}_{f_0}.$$ Let $\xi_t:\mathbb T^2\to\mathbb T^2$ -- be a smooth isotopy such that $\xi_0=id$ and $\xi_1=\xi$. Then the arc $\xi_t f \xi_t^{-1}:\mathbb T^2\to\mathbb T^2$ connects the diffeomorphism $f$ with the diffeomorphism $f_1=\xi f \xi^{-1}\in G$ such that $$c^{u1}_{f_1}=c^{u1}_{f_0},\,c^{u2}_{f_1}=c^{u2}_{f_0}.$$ 

Since the diffeomorphisms $f_1$ and $f_0$ are topologically conjugate on the closures of unstable saddle manifolds, according to the lemma \ref{st2}, there exists an arc without biffurcations connecting $f_1$ with the diffeomorphism $f_2\in G$, which coincides with $f_0$ in some neighborhoods $K_{1}^u,\,K^u_2$ of curves  $c^{u1}_{f_0},\,c^{u2}_{f_0}$.

The set $A=c^{u1}_{f_0}\cup c^{u2}_{f_0}$ an attractor of the diffeomorphisms $f_2$, $f_0$ and the set $U_A=K_{1}^u\cup K^u_2$ is its neighborhood. By construction, the space $\mathbb T^2\setminus U_A$ is homeomorphic to a two-dimensional disk and contains in its interior the points $\alpha$ and $\alpha_{f_2}$. According to \cite [Theorem 3.2, chapter 8]{Hirsh}, there exists a diffeomorphism smoothly isotopic to the identity $\eta:\mathbb T^2\to\mathbb T^2$ such that $$\eta(\alpha_{f})=\alpha,\,\eta|_{U_A}=id.$$ According to the isotopy extension theorem (see, for example, \cite[Theorem 5.8]{Mil}), there exists an arc $\eta_t:\mathbb T^2\to\mathbb T^2$  such that $\eta_0=id$, $\eta_1=\eta$ and $\eta_t|_{U_A}=id$. Then the arc $\eta_t f_2 \eta_t^{-1}:\mathbb T^2\to\mathbb T^2$ connects the diffeomorphism $f_2$ with the diffeomorphism $f_3=\eta f _2\eta^{-1}\in G$ such that $f_2|_{U_A\cup\alpha}=f_0|_{U_A\cup\alpha}$. Moreover, according to Lemma \ref{PR4.3}, we can assume that the diffeomorphism $f_3$ coincides with the diffeomorphism $f_0$ in a neighborhood of the source $\alpha$. Moreover, according to Lemma \ref{statr}, exists an arc without biffurcations connecting the diffeomorphism $f_3$ with the diffeomorphism $f_0$.

2) For a diffeomorphism $f\in G$  such that $J_f\neq E$ let $h=\widehat J^{-1}_f f \widehat J_f:\mathbb T^2\to\mathbb T^2$. Then the diffeomorphism $h$ belongs to the class $G$ and $J_h=E$. According to item 1)  there exists an arc without bifurcations $\zeta_t:\mathbb T^2\to\mathbb T^2$ such that $\zeta_0=h$ and $\zeta_1=f_0$. Then $$H_{f,t}=\widehat J_f \zeta_t \widehat J^{-1}_f:\mathbb T^2\to\mathbb T^2$$ -- is the required isotopy connecting the diffeomorphism $H_{f,0}=f$ with diffeomorphism $H_{f,1}=\widehat J_f f_0 \widehat J^{-1}_f=f_{J_f}$. 
\end{demo}

\subsection{Reduction of a structurally stable diffeomorphism to a linear diffeomorphism in neighborhoods of hyperbolic periodic points}

Let $p$ be a hyperbolic fixed point of a diffeomorphism $f:M^n\to M^n$. {\it A type of point} $p$ is the set of parameters $(q_p,\nu_p,\mu_p)$, where $q_p=\dim~W^u_{p}$, $\nu_p=+1~(-1)$, if $f\vert_{W^u_{p}}$ preserves (reverses) orientation and $\mu_p=+1~(-1)$, if $f\vert_{W^s_{p}}$ preserves (reverses) orientation. According to \cite[Theorem 5.5]{PaMe}, the diffeomorphism $f$ in some neighborhood of a point $p$ of type $(q_p,\nu_p,\mu_p)$ is topologically conjugate to a linear diffeomorphism of the space $\mathbb{R}^n$, defined by
$$A_p=\begin{pmatrix} \nu_p\cdot 2&0&\dots&0&0&0&\dots&0\cr
0&2&\dots&0&0&0&\dots&0\cr & &\ddots\cr
0&0&\dots&2&0&0&\dots&0\cr
0&0&\dots&0&\mu_p\cdot 1/2&0&\dots&0\cr
0&0&\dots&0&0&1/2&\dots&0\cr & & & & & &\ddots\cr
0&0&\dots&0&0&0&\dots&1/2\cr
\end{pmatrix},$$
the number of the rows of $A_p$, containing $2$ (including $\nu_p\cdot 2$), equals $q_{p}$. Denote by $\bar A_p:\mathbb R^n\to\mathbb R^n$ a linear diffeomorphism defined by $A_p$. Let $\mathbb R^u=Ox_1\dots x_{q_p},\,\mathbb R^s=Ox_{q_p+1}\dots x_{n}$, $\bar A^u_p=\bar A_p|_{\mathbb R^u}$ and $\bar A^s_p=\bar A_p|_{\mathbb R^s}$. Then in local coordinates $x^u=(x_1,\dots, x_{q_p})\in \mathbb R^u,\,x^s=(x_{q_p+1},\dots, x_{n})\in \mathbb R^s$ diffeomorphism $\bar A_p$ has a form $$\bar A_p(x^u,x^s)=(\bar A^u_p(x^u),\bar A^s_p(x^s)).$$

\begin{lemma} Let a structurally stable diffeomorphism $\varphi_0: M^n\rightarrow M^n$ has an isolated hyperbolic fixed point $p$, let $(U_0,\psi_0)$ be a local manifold chart of $M^n$ such that $p\in U_0$, $\psi_0(p)=O$ and $U_0$ does not contain non-wandering points of the diffeomorphism $\varphi_0$ other than $p$. Therefore, there exists neighborhoods $U_1, U_2$ of the point $p$, $U_2\subset U_1\subset U_0$ and the arc $\varphi_t:M^n\to M^n,t\in[0,1]$ without bifurcations such that:
 
1) the diffeomorphism $\varphi_t,\,t\in[0,1]$ coincides with the diffeomorphism $\varphi_0$ outside the set $U_{1}$;

2) $p$ is an isolated hyperbolic point for each $\varphi_t$;

3) $W^s_{p}(\varphi_t)=W^s_{p}(\varphi_0)$ and $W^u_{p}(\varphi_t)=W^u_{p}(\varphi_0)$ outside the set $U_{1}$; 

4) the diffeomorphism $\psi_0\varphi_1\psi_0^{-1}$ coincides with the diffeomorphism $\bar A_p$ on the set $\psi_0(U_{2})$. \label{PR4.3}
\end{lemma}
\begin{demo} For $r>0$ let $B_r=\{(x_1,\dots,x_n)\in\mathbb R^n: \sum\limits_{i=1}^nx_i^2\leqslant r^2\}$, $B^u_r=\{(x_1,\dots, x_{q_p})\in\mathbb R^{u}:\sum\limits_{i=1}^{q_p}x_i^2\leqslant r^2\}$ and $B^s_r=\{(x_{q_p+1},\dots, x_{n})\in\mathbb R^{s}:\sum\limits_{i=q_p+1}^nx_i^2\leqslant r^2\}$. 

By virtue of the structural stability of the diffeomorphism $\varphi_0$ any diffeomorphism sufficiently close to $\varphi_0$ in the $C^1$-topology can be joined with $\varphi_0$ by an arc without bifurcations. According to Franks' lemma \cite{Fr71} \footnote{In Franks' lemma, in a neighborhood $U_p$ of the fixed point $p$ of the diffeomorphism $f:M^n\to M^n$ we consider the local chart $(U_p,\psi_p)$ where $\psi^{-1}_p=exp:T_xM^n\to U_p$-- exponential map. Then in these local coordinates the diffeomorphism $f$ has the form $\hat f=exp^{-1}f exp:\mathbb R^n\to\mathbb R^n$. The Franks lemma is that in any neighborhood of a diffeomorphism $f$ there exists a diffeomorphism $g$ having a fixed point $p$ and a linear local representation $\hat g=exp^{-1}g exp$ if it is close enough to $Df_p$. Thus, in any neighborhood of the diffeomorphism $f$ there exists a diffeomorphism $g$,  having a fixed point $p$ and a linear local representation given by a matrix all of whose eigenvalues are pairwise different.} we can assume that the diffeomorphism $\bar\varphi_0=\psi_0\varphi_0\psi_0^{-1}$ in some ball $B_{r_0}\subset\psi_0(U_0)$ coincides with the linear diffeomorphism $\bar\Phi_p:\mathbb R^n\to\mathbb R^n$ given by a matrix $\Phi_p$ all of whose eigenvalues are pairwise different. 
Then the diffeomorphism $\bar\Phi_p$ is smoothly conjugate to the linear diffeomorphism $\bar Q_p$ given by the normal Jordan form $Q_p$ of the matrix $\Phi_p$ (see, for example, \cite[Chapter 3]{gel}). That is, there exists an orientation-preserving diffeomorphism $\xi:\mathbb R^n\to\mathbb R^n$ such that $\bar Q_p=\xi\bar\Phi_p\xi^{-1}$. According to \cite[section 6, Lemma 2]{Mil}, $\xi$ is isotopic to the identity, which means there is an isotopy $\xi_t$ from $\xi_0=id$ to $\xi_1=\xi$.

According to the isotopy extension theorem (see, for example, \cite[Theorem 5.8]{Mil}), there is an isotopy $\Xi_t:\mathbb R^2\to\mathbb R^2$ between the identity map $\Xi_0=id$ and the diffeomorphism $\Xi_t$ coincides with $\xi_t$ on $B_{r_2}$ and is the identity map outside $B_{r_1}$ for some $r_2<r_1< r_0$. 

Thus, the arc $\bar\eta_t =\Xi_t\bar\Phi_p\Xi_t^{-1}:\mathbb R^n\to\mathbb R^n$ connects the diffeomorphism $\bar\eta_0=\bar\Phi_p$ with a diffeomorphism $\bar\eta_1$, coinciding with $\bar Q_p$ on $B_{r_2}$ and with $\bar\Phi_p$ outside $B_{r_1}$. Additionally, $\bar\eta_t$ is an arc without bifurcations $O$, an isolated hyperbolic point for each $\bar\eta_t$ and $W^s_{O}(\bar\eta_t)=W^s_{O}(\bar\Phi_p)$, $W^u_{O}(\bar\eta_t)=W^u_{O}(\bar\Phi_p)$ outside the set $B_{r_1}$.

If $Q_p=A_p$, then the lemma is proved. Otherwise, due to the fact that the eigenvalues of the matrix $Q_p$ are pairwise different, it has a quasi-diagonal form with blocks consisting either of eigenvalues or of matrices of the form     
$\left(\begin{array}{cccc}  
           \alpha & \beta \\  
           -\beta & \alpha \\  
        \end{array}\right),$ where $0<\alpha^2+\beta^2<1$ or $\alpha^2+\beta^2>1$. Then the diffeomorphism $\bar Q_p$ has the form $$\bar Q_p(x^u,x^s)=(\bar Q^u_p(x^u),\bar Q^s_p(x^s)),$$ where $(\bar Q^u_p)^{-1}(B^u_r)\subset int\,B^u_r$ for every disk $B^u_r$ and $\bar Q^s_p(B^s_r)\subset int\,B^s_r$ for every disk $B^s_r$. Choose $r_3<r_2$ in such a way that $B^u_{r_3}\times (\bar Q^s_p)^{-1}(B^s_{r_3})\subset int\,B_{r_2}$. Choose $r^u_4,r^s_4\in(r_3/2,r_3)$ in such a way that $(\bar Q^u_p)^{-1}(B^u_{r_3})\subset int\,B^u_{r_4^u}$ and $\bar Q^s_p(B^s_{r_3})\subset int\,B_{r_4^s}$.

In the proof of the proposition 5.4 of monography\cite{PaMe} arcs $\bar\tau^u_t:\mathbb R^u\to\mathbb R^u,\,\bar\tau^s_t:\mathbb R^s\to\mathbb R^s$ are constructed composing by linear hyperbolic contractions such that
\begin{itemize}
\item $(\bar\tau^u_t)^{-1}(B^u_r)\subset int\,B^u_r$ for any disk $B^u_r$ and $\bar\tau^s_t(B^s_r)\subset int\,B^s_r$ for any disk $B^s_r$;
\item $\bar\tau^u_0=\bar Q^u_p,\,\bar\tau^u_1=\bar A^u_p$ and $\bar\tau^s_0=\bar Q^s_p,\,\bar\tau^s_1=\bar A^s_p$.
\end{itemize}

Consider isotopies $\bar\lambda^u_t=\bar Q^{u}_p\bar\tau^u_t,\,\bar\lambda^s_t=(\bar Q^{s}_p)^{-1}\bar\tau^s_t$, which takes the identity maps $\bar\lambda^u_0=\bar\lambda^s_0=id$ to diffeomorphism $\bar\lambda^u_1=\bar Q^u_p(\bar A^u_p)^{-1},\,\bar\lambda^s_1=(\bar Q^s_p)^{-1}\bar A^s_p$, respectively. By construction, $\bar\lambda^u_t (B^u_{r_3})\subset \bar Q^u_p(B_{r^u_4})$ and $\bar\lambda^s_t (B^s_{r_3})\subset (\bar Q^s_p)^{-1}(B_{r^s_4})$ for each $t\in[0,1]$. Thus, by virtue of the isotopy extension theorem (see e.g., \cite[Theorem 5.8]{Mil}), there exist isotopies $\bar\Lambda^u_t:\mathbb R^u\to\mathbb R^u,\,\bar\Lambda^s_t:\mathbb R^s\to\mathbb R^s$ taking the identity map
$\bar\Lambda^u_0=\bar\Lambda^s_0=id$ which coincide with $\bar\lambda^u_t,\,\bar\lambda^s_t$ on $B^u_{r_3},\,B^s_{r_3}$ and are exactly identity outside $\bar Q^u_p(B_{r^u_4}),\,(\bar Q^s_p)^{-1}(B_{r^s_4})$, respectively. Let $$\bar\Lambda_t(x^u,x^s)=((\bar\Lambda^u_t)^{-1}\bar Q^u_p(x^u),\bar Q^s_p\bar\Lambda^s_t(x^s)).$$
Let us denote an arc coinciding with $\bar\eta_1$ outside $B_{r_2}$ and with $\bar\Lambda_t$ on $B_{r_2}$ by $\bar\zeta_t:\mathbb R^n\to\mathbb R^n$. Choose $r_5<r_4$, such that $B_{r_5}\subset\bar Q^u_p(B^u_{r_3})\times B^s_{r_3}$. Let $\bar U_2=B_{r_5}$, $\bar U_1=B_{r_1}$ and $$\bar\varphi_t=\bar\eta_t*\bar\zeta_t.$$
Then $\bar\varphi_t$ is an arc without bifurcations, which coincide with $\bar\varphi_0$ outside $\bar U_{1}$, $O$ is an isolated hyperbolic point for each $\bar\varphi_t$, $W^s_{O}(\bar\varphi_t)=W^s_{O}(\bar\varphi_0)$, $W^u_{O}(\bar\varphi_t)=W^u_{O}(\bar\varphi_0)$ outside the set $\bar U_{1}$ and the diffeomorphism $\bar\varphi_1$ coinsedes with the diffeomorphism $\bar A_p$ on $\bar U_2$. So, the arc $\varphi_t:M^n\to M^n$ which coincides with $\psi_0^{-1}\bar\varphi_1\psi_0$ on $U_0$ and with $\varphi_0$ outside $U_0$ satisfies all the conditions of the lemma in $U_1=\psi_0^{-1}(\bar U_1)$ and $U_2=\psi_0^{-1}(\bar U_2)$.
\end{demo}

\begin{collary} The statement analogous to \ref{PR4.3} holds, if point $p$ is periodic with period $m$. To prove that it is sufficient to consider isotopy for the diffeomorphism $f^m$ with fixed point $p$, which exists according to lemma \ref{PR4.3}, and extend it along the neighbourhood of point $p$ orbit with $f$. \label{cPR4.3}
\end{collary}

\section{Straightening of invariant saddle manifolds}

The invariant manifold of the saddle point of a Morse-Smale diffeomorphism is always a smooth submanifold. However, its closure, even in a gradient-like case, may not even be a topological submanifold (Fig. \ref{42}). In this section, we present statements and facts that allow us to connect any gradient-like 2-diffeomorphism by an arc without bifurcations with a diffeomorphism for which the closures of the invariant manifolds of all saddle points are smooth submanifolds.

\begin{proposition} \label{Rol} (\cite[Theorem 13, p. 25]{Rolf} and \cite{Mu}, \cite[Proposition 10.60]{GrPo} and \cite[Statement 10.51]{GrPo}) 
Suppose that on the torus $\mathbb T^2$ there are $p\in\mathbb N\cup \{0\}$ pairwise disjoint closed smooth curves $c_1,\dots,c_{p}$, of homotopy type $<1,0>$ and $q\in\mathbb N\cup \{0\}$ pairwise disjoint closed smooth curves $d_1,\dots,d_{q}$, of homotopy type $<0,1>$ such that each pair of curves $c_i,\,d_j$ intersects transversally at one point. Let $C_p$ be a disjunct union of $p$ circles of the form $\mathbb S^1\times \{w\},\,w\in\mathbb S^1$ and $D_q$ be the disjoint union of $q$ circles of the form $\{z\}\times \mathbb S^1,\,z\in\mathbb S^1$.

Then there is an isotopic to identity diffeomorphism $\xi:\mathbb T^2\to\mathbb T^2$ such that $\xi(c_1\cup\dots\cup c_p)=C_p$ and $\xi(d_1\cup\dots\cup d_q)=D_q$.
\end{proposition}

\begin{lemma}\label{L5} Let a diffeomorphism $\varphi_0:M^2\to M^2$ has a hyperbolic sink $\omega_0$ and hyperbolic saddles $\sigma_1,\dots,\sigma_k$ such that their unstable separatrices $\gamma^1_{\varphi_0},\dots,\gamma^k_{\varphi_0}$ lies in the sink basin $W^s_{\omega_0}$ and has the same period $m$ as the sink $\omega_0$. Let $(U_0,\psi_0)$ be a local chart of the manifold $M^2$ such that $\omega_0\in U_0$,  $\psi_0(\omega_0)=O$ and $\varphi_0^m(U_0) \subset U_0$. Let $L_k\subset\mathbb{R}^2$ be a sheaf of rays $l_1,\dots,l_k$, which in polar coordinates $(\rho,\theta)$ has the form $l_i=\{(\rho,\theta)\in\mathbb R^2:\theta=\theta_i\}$, $\theta_i\in[0,2\pi)$. Then there are neighborhoods $V_1, V_2$ of the point $\omega_0$ such that $V_2\subset V_1\subset U_0$ and the arc $\varphi_t:M^2\to M^2,\,t\in[0,1]$ without bifurcations with the following properties:

1) the diffeomorphism $\varphi_t,\,t\in[0,1]$ coincides with the diffeomorphism $\varphi_0$ outside the set $\bigcup\limits_{k=0}^{m-1}\varphi_0^k(V_{1})$ and $\bigcup\limits_{k=0}^{m-1}\varphi_0^k(\omega_0)$ is a hyperbolic sink orbit of period $m$ for all $\varphi_t$;

2) $\psi_0(\bigcup\limits_{i=0}^{k} \gamma^i_{\varphi_1})\cap V_2)\subset L_k$, where $\gamma^i_{\varphi_1}$ are unstable separatrices of saddles $\sigma_i$ with respect to the diffeomorphism $\varphi_1$.
\end{lemma}
\begin{demo}
\begin{figure}[h]
\centerline{\includegraphics[width=12 cm]{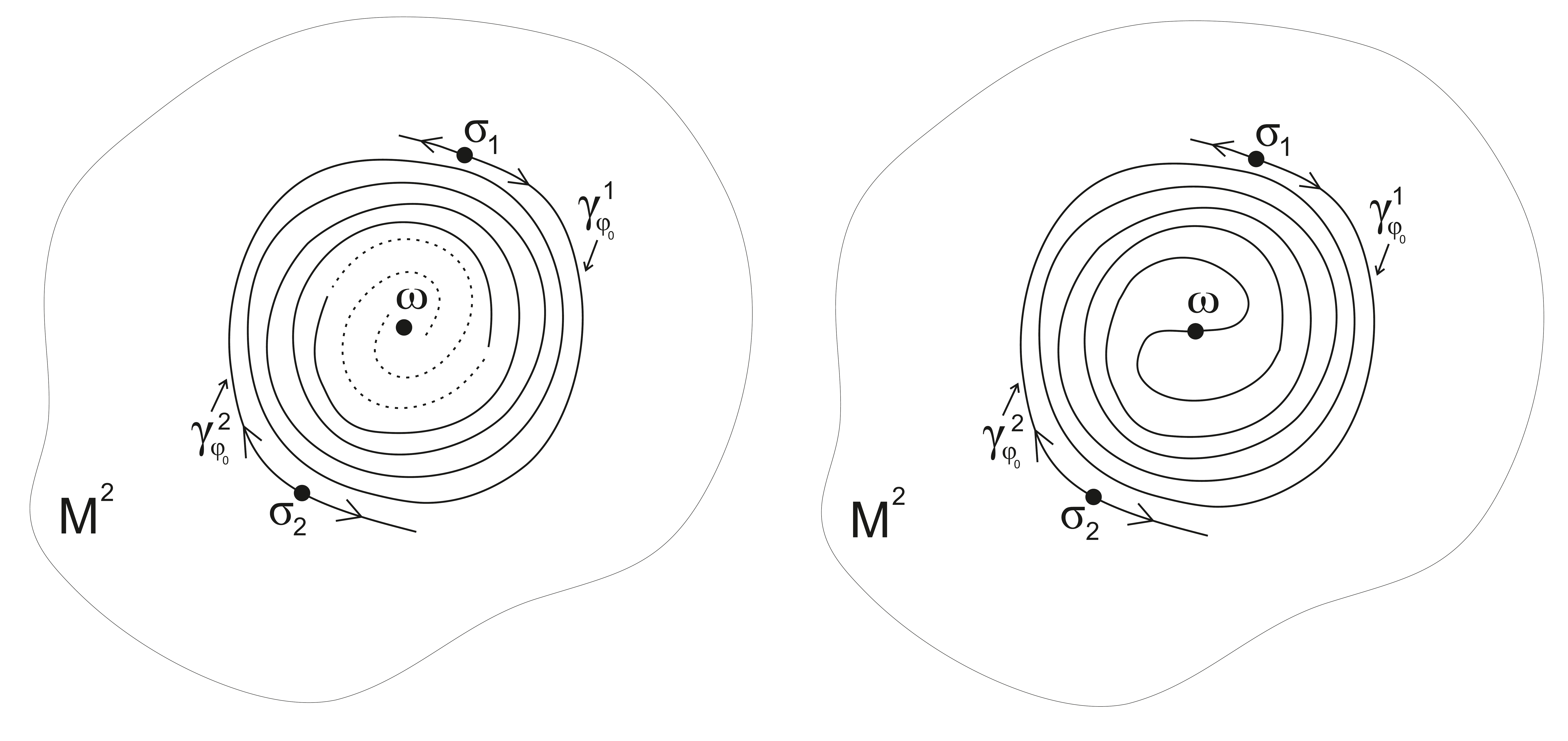}}
\caption{\small Straightening the separatrix}\label{42}
\end{figure}
Let $\phi_0=\varphi_0^m$ and $\bar\phi_0=\psi_0\phi_0\psi_0^{-1}:\mathbb R^2\to\mathbb R^2$. 
Denote by $O(0,0)$ the origin of coordinates on the plane $\mathbb R^2$. For any $r>0$ we put $B_r=\{(x,y)\in\mathbb R^2:x^2+y^2<r^2\}$. According to lemma \ref{PR4.3}, without loss of generality, we can assume that $\bar\phi_0=g,$  where $ g(x,y)=\left(\frac x2, \frac y2\right)$ on disk $B_{2r_0}$ for some $r_0>0$. Let $K_0=B_{2r_0}\setminus B_{r_0}$,$\gamma^i_{\bar\phi_0}=\psi_0(\gamma^i_{\varphi_0})$.

Denote by $E_g$ the set of contractions $\bar\phi:\mathbb R^2\to\mathbb R^2$, coinciding with $\bar\phi_0$ outside $B_{2r_0}$ and with $g$ with $B_{r_{\bar\phi}}$, where $r_{\bar\phi}\leq 2 r_0$. For any $\bar\phi\in E_g$ let $$\gamma^i_{\bar\phi}=\bigcup\limits_{k\in\mathbb Z}\bar\phi^k(\gamma^i_{\bar\phi_0}\cap K_0).$$ By construction, the $\bar\phi$-invariant curves $\gamma^i_{\bar\phi}$ coincide with $\bar\phi_0$-invariant curves $\gamma^i_{\bar\phi_0}$ outside the disk $B_{r_0}$. Then, to prove the lemma, it suffices to construct an arc of contractions $\bar\phi_t:\mathbb R^2\to\mathbb R^2,t\in[0,1]$ such that

$\bar 1)$  the diffeomorphism $\bar\phi_t,t\in[0,1]$ coincides with the diffeomorphism $\bar\phi_0$ outside the set $B_{r_0}$;

$\bar 2)$  $\psi_0(\bigcup\limits_{i=0}^{k} \gamma^i_{\varphi_1})\cap V_2)\subset L_k$. 

The arc $\varphi_t:M^2\to M^2$ is obtained from the arc $\bar\phi_t$ as in lemma \ref{PR4.3}, if we put $V_1=\psi^{-1}_0(B_{r_0})$ and $V_2=\psi_0^{-1}(B_{r_{\bar\phi_1}})$.

To construct the arc $\bar\phi_t$ we introduce the following notations for any diffeomorphism $\bar\phi\in E_g$. 

Represent the 2-torus $\mathbb T^2$ as the orbit space of the diffeomorphism $g$ on the set $\mathbb R^2\setminus O$ and
denote by $p:\mathbb R^2\setminus O\to\mathbb T^2$ the natural projection. Let $\hat a=p(OX^+)$ and $\hat b=p(\mathbb S^1)$ be generators of $\mathbb T^2$.
Let $\hat l_i=p(l_i)$ and $K_{\bar\phi}=B_{r_{\bar\phi}}\setminus B_{r_{\bar\phi}/2}$, $\hat\gamma^i_{\bar\phi}=p(\gamma^i_{\bar\phi}\cap K_{\bar\phi})$. Then the curve $\hat\gamma^i_{\bar\phi}$ is the knot on the torus $\mathbb T^2$ with the homotopy type $<1,-n_{\bar\phi}>,\,n_{\bar\phi}\in\mathbb Z$ in the basis $\hat a,\hat b$ (see, for example, \cite{GrPo}).

The arc $\bar\phi_t$ is the smooth product of $\eta_t$ and $\zeta_t$ where

I) the arc $\eta_t,t\in[0,1]$  consists of contractions, it coincides with $\bar\phi_0$ outside the set $B_{r_0}$ and
it joins the diffeomorphism ${\eta}_0=\bar\phi_0$ with some diffeomorphism ${\eta}_1\in E_g$ such that the knot $\hat{\gamma}_{\eta_1}$ and $\hat{\xi}_{\eta_1}$ has the homotopy type $<1,0>$ in the basis $\hat a,\hat b$;

II) the arc $\zeta_t\in E_g,t\in[0,1]$  joins the diffeomorphism $\zeta_0=\eta_1$ to the diffeomorphism $\zeta_1$ such that $\hat\gamma^i_{\zeta_1}=\hat l_i$.  

\begin{figure}[h]
\centerline{\includegraphics[width=17 cm]{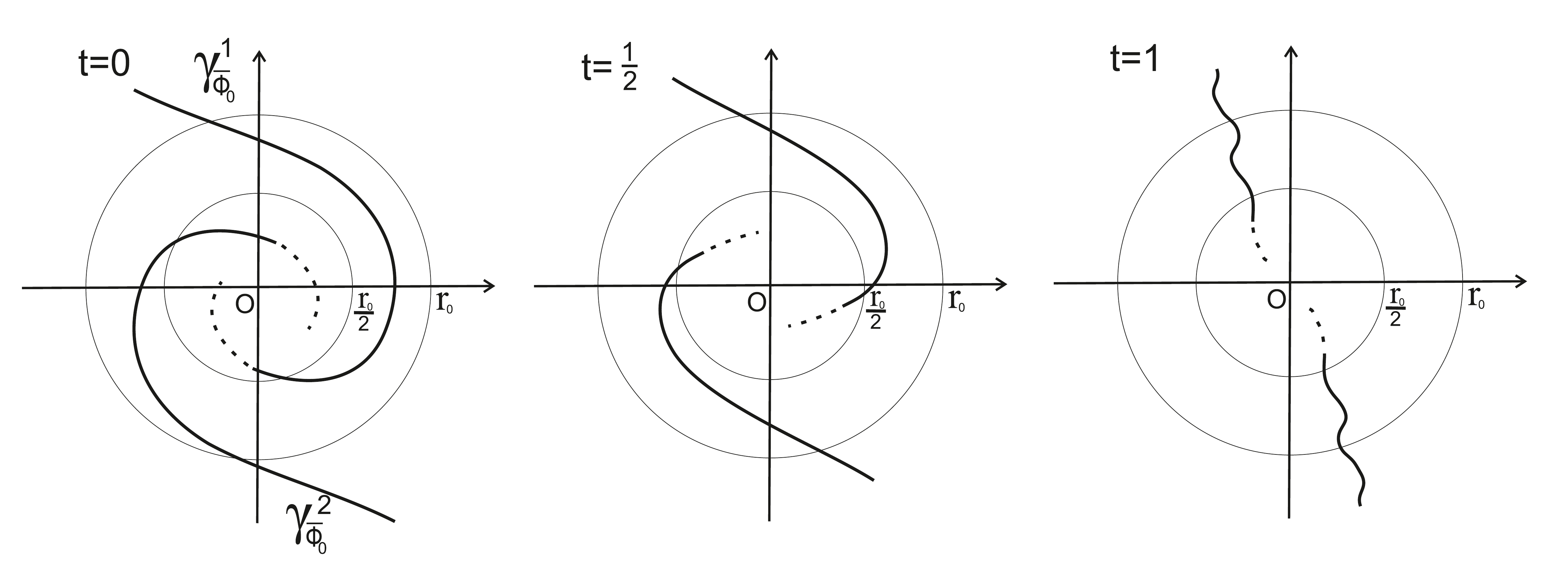}}
\caption{\small Illustration to Lemma \ref{L5}, Part I)}\label{421}
\end{figure}

I) If $n_{\bar\phi}= 0$, then let $\eta_t=\bar\phi_0$ for every $t\in[0,1]$. Otherwise define the diffeomorphism $\theta_t,t\in[0,1]:\mathbb R^2\to\mathbb R^2$ in polar coordinates $\rho,\varphi$ in such a way that $\theta_t(O)=O$ and ${\theta}_t(\rho e^{i\varphi})=\rho e^{i(\varphi+\varphi_t(\rho))}$ where $\varphi_t(\rho)$ is a smooth monotonic function equal to $2n_{\bar\phi}\pi t$ for $\rho\leqslant\frac{1}{2}$ and equal to $0$ for $\rho\geqslant{1}$. 
Then $\eta_t=\theta_t \bar\phi_0:\mathbb R^2\to\mathbb R^2$ is the desired arc (Fig. \ref{421}). 

II) By construction, the diffeomorphism $\eta_1\in E_g$ and the knot $\hat \gamma^i_{\eta_1}$ has the homotopy type $<1,0>$ in the basis $\hat a,\hat b$. Let $\hat L_k=p(L_k)$. According to \ref{Rol}, there is a diffeomorphism $\hat h:\mathbb T^2\to\mathbb T^2$ which is smoothly isotopic to the identity and such that $\hat h(\bigcup\limits_{i=1}^k\hat\gamma^i_{\eta_1})=\hat L_k$. For $r>0$ let $K_r=B_r\setminus B_{r/2}$. Pick an open cover $D=\{D_1,\dots, D_q\}$ of the torus $\mathbb{T}^2$ such that the connected component $\bar D_i$ of the set $p^{-1}(D_i)$  is a subset of $K_{r_i} $ for some $r_i<\frac{r_{i-1}}{2}$ and $r_1\leq r_0/2$. According to \cite[Lemma de fragmentation]{Ba} there are diffeomorphisms $\hat{w}_ {1}, \dots, \hat {w} _ {q}: \mathbb {T}^2\to\mathbb{T}^2$ smoothly isotopic to the
identity, with the following properties:

i) for each $i\in\{1,\dots, q\}$ there exists a smooth isotopy $\{\hat{w}_{i,t}\}$ which is the identity outside $D_{i}$ and which joins the identity and $\hat{w}_{i}$;

ii) $\hat{h}=\hat{w}_{1} \dots \hat {w} _ {q} $.

Let ${w}_{i,t}:\mathbb{R}^2\to\mathbb{R}^2$ be a diffeomorphism which coincides with $(p\vert_{K_{r_i}})^{-1}\hat{w}_{i,t}p$ on $K_{r_i}$ and which coincides with the identity outside $K_{r_i}$ (Fig. \ref{422}). 
Then the desired arc is defined by $${\zeta}_{t}=\nu_1 {w}_{1,t}\dots {w}_{q,t}.$$ 

\begin{figure}[h]
\centerline{\includegraphics[width=17 cm]{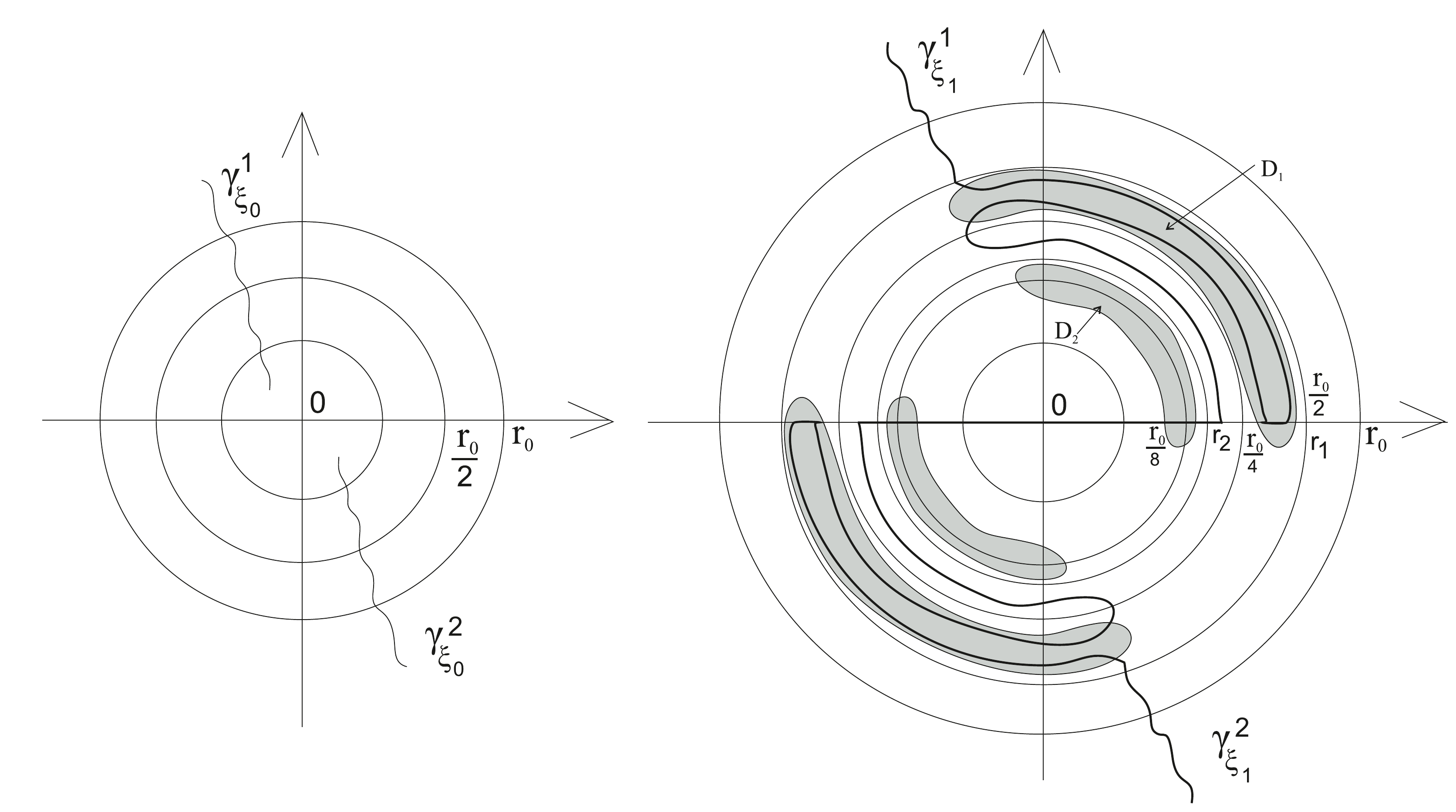}}
\caption{\small Illustration to Lemma \ref{L5}, Part II)}\label{422}
\end{figure}
\end{demo}

\subsection{Changes in dynamics in the vicinity of the attractor}
In this section, we prove that the dynamics of a gradient-like 2-diffeomorphism in a neighborhood of a smoothly embedded closure of an unstable manifold of a saddle point can be replaced by any topologically conjugate one by passing along an arc without bifurcations.

\begin{lemma}\label{st2} Let 
\begin{itemize}
\item $\phi_0:M^2\to M^2$ be a orientation-preserving Morse-Smale diffeomorphism defined on an orientable surface $M^2$, 
\item $\sigma_1,\dots,\sigma_n$ be saddle fixed points of the diffeomorphism $\phi_0$, the curve $c=\bigcup\limits_{i=1}^n cl\,W^u_{\sigma_i}$ be smoothly embedded in $M^2$ and $K_0$ be a basin of the attractor,

\item a diffeomorphism $\phi_1$ be topologically conjugate with the diffeomorphism $\phi_0$ on the curve $c$ and in some of its neighborhood. 
\end{itemize} 
Then there is an arc without bifurcations $\varphi_t:M^2\to M^2,\,t\in[0,1]$ and a neighborhood $K_*\subset K_0$ of the curve $c$ with the following properties:
\begin{enumerate}
\item $\varphi_0=\phi_0$, $\varphi_1$ coincides with $\phi_1$ on $K_{*}$ and coincides with $\phi_0$ outside $K_0$;
\item $\phi_1$ coincides with $\phi_0$ outside $K_0$ and if $\phi_1$ coincides with $\phi_0$ in some neighborhood of the sink on the curve $c$ then $\varphi_t$ coincides with $\phi_0$ in this neighborhood.
\end{enumerate}
\end{lemma}
\begin{demo} Consider the case when the curve $c$ is closed (in the case when the curve $c$ is not closed, the proof is similar). Then there is a smooth embedding $\nu:\mathbb S^1\times[-1,1]\to M^2$ such that $\nu:\mathbb S^1\times\{0\}=c$ and $K_0=\nu(\mathbb S^1\times[-1,1])$. 

Then the desired arc $\varphi_t$ will be the smooth product of arcs without bifurcations $\mu_t$ and $\delta_t$ where

I) the arc $\mu_t,t\in[0,1]$ joins the diffeomorphism $\mu_0=\varphi_0$ with the diffeomorphism $\mu_1$ such that $\Omega_{\eta_1|_c}=\Omega_{\phi_1|_c}$;

II) the arc $\delta_t,t\in[0,1]$ joins the diffeomorphism $\delta_0=\mu_1$ with the desired diffeomorphism $\delta_1=\varphi_1$.

I) If $\Omega_{\phi_0|_c}=\Omega_{\phi_1|_c}$ then $\mu_t=\phi_0$ for all $t\in[0,1]$. Otherwise, by the condition the diffeomorphisms $w_0={\phi_0}|_c,\,w_1={\phi_1}|_c$ are topologically conjugate to rough transformations of the circle. According to \cite{Ma},  they have the same number of hyperbolic attracting and repelling periodic points alternating on a circle. Then there are three possible cases where these points are located on the circle $c$: 
1) there are periodic points $p_0^0,p_0^1$ of diffeomorphisms $w_0,\,w_1$, respectively, which coincide and have the same stability; 
2) there exists an arc $\beta\subset c$ such that $int\,\beta$ contains periodic points $p_0^0,p_0^1$ of diffeomorphisms $w_0,\,w_1$, respectively, of the same stability and does not contain other periodic points of diffeomorphisms $w_0,\,w_1$; 
3) there exists an arc $\beta\subset c$ such that $int\,\beta$ contains all periodic points of both diffeomorphisms $p_0^0,p_1^0,p_0^1,p_1^1,\dots,p_{2n-2}^0,p_{2n-1}^0,p_{2n-2}^1,p_{2n-1}^1$ located on the arc in the specified order, where $p_{2i}^0,p_{2i}^1$ are attracting points, $p_{2i-1}^0,p_{2i-1}^1$ are repelling points of diffeomorphisms $w_0,\,w_1$, respectively, for $i\in\{0,\dots,n-1\}$. 

In case 1) let $\beta=c\setminus U_0$, where $U_0$ is the arc of the circle $c$ containing the point $p_0^0$ and not containing any other periodic points of the diffeomorphisms $w_0,\,w_1$. Choose a tubular neighborhood $N(\beta)\subset K_0$ of the arc $\beta$ so that there exists a diffeomorphism $h: N(\beta)\to\Pi$ where $\Pi= [0,1]\times[-1,1]$. We renumber $p_1^0,p_2^0,\dots,p_{2n-1}^0;\,p_1^1,p_2^1,\dots,p_{2n-1}^1$ periodic points of diffeomorphisms $w_0,\,w_1$, respectively, on the arc $\beta$ so that $0<h(p_1^0)<h(p_2^0)<\dots<h(p_{2n-1}^0)<1;\,0<h(p_1^1)<h(p_2^1)<\dots<h(p_{2n-1}^1)<1$. 
\begin{figure}[h]
\centerline{\includegraphics[width=6 cm]{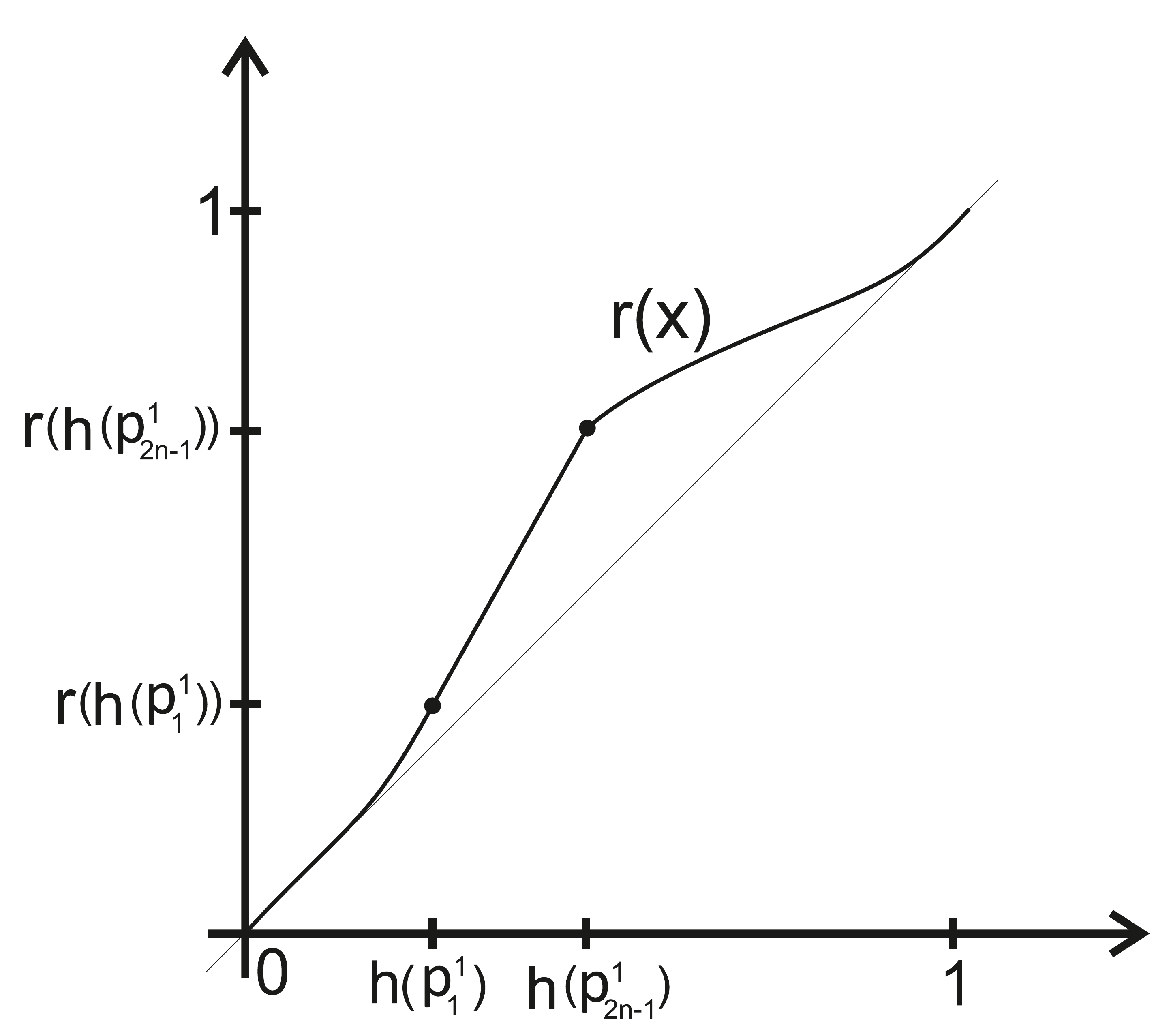}}
\caption{\small Graph of the map $r(x)$}\label{q(x)}
\end{figure}

Consider a smooth monotonic increasing function $r(x):[0,1]\to[0,1]$ which is identical in some neighborhoods of the points $0,\,1$ and on the set $h(\nu(B))$ such that $r(h(p_i^0))=h(p_i^1)$  (Fig. \ref{q(x)}). Let $r_t(x)=tr(x)+(1-t)x,\, x\in[0,1],\,t \in [0,1]$ and $r_{t,y}(x)=(1-y)r_t(x)+yx, x\in [0,1],\,t\in[0,1],\,y \in [0,1]$. We define a smooth isotopy $R_t:\Pi\to\Pi$ by the formula $$R_t(x,y)=(r_{t,y^2}(x),y).$$
Since the isotopy $R_t$ is identical on $\partial\Pi$ for all $t\in[0,1]$, there exists a smooth isotopy $\rho_t:M^2\to M^2$ which is identical outside $N(\beta)$ and coincides with $h^{-1}R_th$ on $N(\beta)$. Then $\mu_t=\rho_tf\rho^{-1}_t$ is the required isotopy that superpose the periodic points of the diffeomorphism $w_0$ with the periodic points of the diffeomorphism $w_1$ on the arc $\beta$.

In case 2) using the technique above, the desired isotopy $\mu_t$ is composed of an isotopy that superpose the point    $p_0^1$ on the arc $\beta$ and an isotopy that combines the periodic points of the diffeomorphism $w_0$ with periodic points of the diffeomorphism $w_1$ on the complementary arc to $\beta$.

In case 3) the desired isotopy $\mu_t$ superpose the periodic points of the diffeomorphism $w_0$ with the periodic points of the diffeomorphism $w_1$ while preserving the order on the arc $\beta$.

If $\phi_1$ coincides with $\phi_0$  in some neighborhood of the sink on the curve $c$, then, evidently, the arc $\mu_t$ can be constructed to coincide with $\phi_0$ in this neighborhood.

II) If ${\mu_1}$ coincides with $\phi_1$ in some neighborhood of the circle $c$, then $\delta_t=\mu_1$ for all $t\in[0,1]$. Otherwise, since the circle $c$ is an attractor for the diffeomorphism $\phi_1$, then there exists a smooth annulus $K\subset int\,K_{0}$, which is the trapping neighborhood of the attractor $c$ of the diffeomorphism $\phi_1$. Choose $0<\varepsilon_*<1$ so that for the annulus $K_*=\nu(\mathbb S^1\times[-\varepsilon_*,\varepsilon_*])$ satisfies the inclusion $K_{*}\subset int\,(\phi^2_1(K)\cap \mu^2_1(K_0))$. From the annulus theorem (see, for example, \cite{Rolf}),  there exists the diffeomorphism $\gamma:K\to K_0$ such that $\gamma|_{K_{\varepsilon_*}}=id$ and $\gamma(\phi^i_1(K))=\mu^i_1(K_0),\,i=1,2$. Define on $\mathbb S^1\times[-1,1]$ the diffeomorphisms $\tilde\phi_1=\nu^{-1}\gamma\phi_1\gamma^{-1}\nu$ and $\tilde\mu_1=\nu^{-1}\mu_1\nu$. Let $$\tilde\xi_t=(1-t)\tilde\mu_1+t\tilde\phi_1.$$ Then the arc $\tilde q_t=\tilde\mu_1^{-1}\tilde\xi_t$ joins the identity map $\tilde q_0=id$ and the diffeomorphism $\tilde q_1=\tilde\mu_1^{-1}\tilde\phi_1$ and $\tilde q_t(\tilde\mu_1(\mathbb S^1\times[-1,1]))\subset \mathbb S^1\times(-1,1)$. According to  isotopy extension theorem (see, for example, \cite[Theorem 5.8]{Mil}), there exists an isotopy $\tilde Q_t:\mathbb S^1\times[-1,1]\to\mathbb S^1\times[-1,1]$, coincides with $\tilde q_t$ on $\tilde\mu_1(\mathbb S^1\times[-1,1])$ and is identical on $\mathbb S^1\times\{-1,1\}$. The desired isotopy which $\delta_t:M^2\to M^2$, coincides with $\mu_1$ outside $K_0$ and coincides with $\mu_1\nu{\tilde Q}_t\nu^{-1}$ on $K_0$.

If $\phi_1$ coincides with $\mu_1$ in some neighborhood of the sink on the curve $c$, then evidently the annulus $K$ can be chosen to coincide with the annulus $K_0$ and the diffeomorphism $\gamma$ is identical in this neighborhood.
\end{demo}

\subsection{Changing dynamics in a wandering set}
The closures of unstable saddle manifolds of any gradient-like 2-diffeomorphism form a connected attractor of this diffeomorphism, and the repeller dual to it consists of all sources (see, for example, \cite{GMPZ2010}). In this section, we will prove that any two diffeomorphisms that coincide in some neighborhoods of such an attractor and a repeller can be connected by an arc without bifurcations.

\begin{lemma}\label{statr} Let the gradient-like 2-diffeomorphisms $\varphi_0,\,\varphi_1:M^2\to M^2$ coincide on the closures of unstable saddle manifolds and in some of their neighborhoods, as well as in neighborhoods of sources. Then there is a bifurcation free arc connecting them.
\end{lemma}
\begin{demo} Let $A$ be an attractor of diffeomorphisms $\varphi_0,\,\varphi_1$,  formed by the closures of the unstable manifolds of all saddle points, and $R$ is its dual repeller consisting of sources. Let $K_A$ be a trapping neighborhood of attractor $A$. Then each connected component $M^2\setminus K_A$ is a two-dimensional disk lying in the basin of some source from the repeller $R$. We pick in the repeller the sources $\alpha_1,\dots,\alpha_k$ with pairwise disjoint orbits of periods $m_1,\dots, m_k$, respectively, so that $R=\bigcup\limits_{i=1}^k\bigcup\limits_{j=0}^{m_i-1}\varphi_0^j(\alpha_i)$. 

We construct an arc $H_{\varphi_0,\alpha_1,t}$ without bifurcations with the following properties:

- $H_{\varphi_0,\alpha_1,t}$ connects the diffeomorphism $\varphi_0$ with the diffeomorphism $\varphi_{\alpha_1}=H_{\varphi_0,\alpha_1,1}$ such that $\varphi_{\alpha_1}=\varphi_1$ on $K_A\cup\bigcup\limits_{j=0}^{m_i-1}\varphi_0^j(W^u_{\alpha_1})$;

- $H_{\varphi_0,\alpha_1,t}=\varphi_0$ outside $K_A\cup\bigcup\limits_{j=0}^{m_i-1}\varphi_0^j(W^u_{\alpha_1})$.

Analogously one constructs the arc $H_{\varphi_{\alpha_1},\alpha_2,t},\dots, H_{\varphi_{\alpha_{k-1}},\alpha_k,t}$. 

Then the desired arc $\varphi_t$ will be the smooth product of these arcs, that is $$\varphi_t=H_{\varphi_0,\alpha_1,t}*\dots*H_{\varphi_{\alpha_{k-1}},\alpha_k,t}$$ 

{\bf Construction of the arc $H_{\varphi_0,\alpha_1,t}$.} For definiteness, we assume that the source $\alpha_1$ is fixed.

In the case when the point $\alpha_1$ is periodic with period $m_1$ suffices to construct an isotopy for the diffeomorphism $f^{m_1}$ with a fixed point $\alpha_1$  and to extend it along the neighborhood of the orbit of the point $\alpha_1$ with diffeomorphism $\varphi_0$ 

Denote by $\Phi_{\varphi_1}$ the set of gradient-like diffeomorphisms of the surface  $M^2$ that coincide with $\varphi_1$ on $K_A$ and in some neighborhood of the point $\alpha_1$ (in particular, the diffeomorphism $\varphi_0$ belongs to the set $\Phi_{\varphi_1}$). Then any diffeomorphism $\varphi\in\Phi_{\varphi_1}$ is smoothly conjugate on $W_{\alpha_1}^u\setminus\alpha_1$ with $\varphi_1$ by the diffeomorphism $$\rho_{\varphi}(x)=\varphi_1^{-k} (\varphi^k(x)),$$ where $k\in\mathbb Z$ such that $\varphi^k(x)\in K_A$ for $x \in W_{\alpha_1}^u\setminus\alpha$. Thus $\varphi_1=\rho_{\varphi} \varphi_0 \rho_{\varphi}^{-1}$ on $W_{\alpha_1}^u\setminus\alpha_1$.

For the diffeomorphism $\varphi$ сonsider two subcases: I) $\rho_{\varphi}$ be smoothly extended to the point $\alpha_1$ by the condition $\rho_{\varphi}(\alpha_1)=\alpha_1$, II) $\rho_{\varphi}$ does not smoothly extended to the point $\alpha_1$.

In case I) according to \cite{Sm-S2} the diffeomorphism $\rho_{\varphi}$ is isotopic to the identity. 

Moreover, according to Thom’s isotopy extension theorem, on $W_{\alpha_1}^u$ there exists a smooth isotopy $\varrho_{\varphi,t}$ such that $\varrho_{\varphi,0}=id,\,\varrho_{\varphi,1}=\rho_{\varphi}$ and $\varrho_{\varphi,t}=id,\,t\in[0,1]$ in some neighborhood $K_\varphi\subset K_A$ of attractor $A$. The isotopy $\varrho_{\varphi,t}$ extends to isotopy of ambient surface $\varrho_{\varphi,t}:M^2\to M^2$ and coincides with identity map outside $W^u_{\alpha_1}$. Then a bifurcation free arc $$H_{\varphi,\alpha_1,t}=\varrho_{\varphi,t}^{-1}\varphi\varrho_{\varphi,t}$$ connects the diffeomorphism $\varphi$ with the diffeomorphism $\varphi_1$.

In case II) let us apply the following construction. Let $(U_{\alpha_1}, \psi_{\alpha_1}), \psi_{\alpha_1} : U_{\alpha_1} \to \mathbb R^2$, $\psi_{\alpha_1}=O$ be a local chart of $M^2$. We consider on $\mathbb R^2$ diffeomorphisms $\bar{\varphi}=\psi_{\alpha_1} \varphi \psi_{\alpha_1}^{-1}$, $\bar{\varphi_1}=\psi_{\alpha_1} \varphi_1 \psi_{\alpha_1}^{-1}$ and $\bar{\rho}_\varphi=\psi_{\alpha_1} \rho_\varphi \psi_{\alpha_1}^{-1}$, which are local representations of the diffeomorphisms $\varphi,\varphi_1$ and $\rho$, respectively. Since the point $O$  is a hyperbolic source of the diffeomorphisms $\varphi$ and $\varphi_1$, there exists a 2-disk $B_\varphi\ni O$ such that $\bar\varphi^{-1}(B_\varphi)\subset int\, B_\varphi$ and the annulus $K_\varphi$ is a fundamental domain of the restriction of the diffeomorphism $\bar \varphi$ to $int\,B_\varphi\setminus\{O\}$.

Represent the 2-torus $\mathbb T^2$ as the orbit space $(int\, B_\varphi\setminus\{O\})/{\bar \varphi}$. Denote by $p_\varphi:B_\varphi\setminus\{O\}\to\mathbb T^2$ the natural projection. Then the curve $b=p_\varphi(\partial B_\varphi)$ has the homotopy type <0,1> and it is possible to uniquely define a curve $a$ of homotopy type <1,0> such that the curves $a,b$ are the generators of the fundamental group $\pi_1(\mathbb T^2)$. 
Since $\bar\rho_\varphi$ maps the orbits $\bar\varphi$ into the orbits $\bar\varphi_1$ and $K_\varphi$ is a common fundamental domain for $\bar\varphi,\bar\varphi_1$ on $int\,\ B_\varphi\setminus\{O\}$,  then $\bar\rho_\varphi$ is projected onto $\mathbb T^2$ by the formula $\hat{\rho}_\varphi=p_\varphi\bar\rho_\varphi p_\varphi^{-1}$. Then the induced isomorphism $\hat{\rho}_{\varphi*}:\pi_1(\mathbb T^2)\to\pi_1(\mathbb T^2)$  preserves the homotopy class of the generator $a$ and, therefore, is given by a matrix $$\left(\begin{array}{cccc} 1 & 0 \\ n_\varphi & 1 \\ \end{array}\right)$$ for some integer $n_\varphi$.

Consider two subcases: a) $n_\varphi=0$; b) $n_\varphi\neq 0$. 

In the case a) let us construct an arc without bifurcations $\nu_{\varphi,t}:M^2\to M^2$ such that $\nu_{\varphi,0}=\varphi$, $\nu_{\varphi,1}\in \Phi_{\varphi_1}$ and the diffeomorphism $\rho_{\nu_{\varphi,1}}$ smoothly continues to the source $\alpha_1$. Then an arc without bifurcations $$H_{\varphi,\alpha_1,t}=\nu_{\varphi,t}*(\varrho_{\nu_{\varphi,1},t}^{-1}\nu_{\varphi,1}\varrho_{\nu_{\varphi,1},t})$$ connects the diffeomorphism $\varphi$ with the diffeomorphism $\varphi_1$. 

Let us describe the construction of the arc  $\nu_{\varphi,t}$.

We choose a cover $U=\{U_1,\dots, U_q\}$ of the torus $\mathbb{T}^2$, consisting of 2-disks such that some connected component of the set $p_\varphi^{-1}(U_i)$ is a subset of $K_i=B_i\setminus \bar\varphi^{-1}(B_i)$, obtained from a 2-disc $B_i$ such that $O\in B_i\subset\bar\varphi^{-1}(B_\varphi),\,i=1,\dots,q$.

According to the lemma of fragmentation \cite{Ba} there exist smoothly isotopic to the identity diffeomorphisms $\hat {w} _ {1}, \dots, \hat {w} _ {q}: \mathbb {T}^2\to\mathbb{T}^2$ such that

i) for each $i =\overline {1, q}$ there exists $U _ {j (i)} \in U$ such that for each $t\in [0,1] $ the map $\hat{w}_{i, t} $ is identical outside $U_{j (i)}$, where $\{\hat{w}_{i,t}\}$ is the smooth isotopy between the identity map and $\hat{w}_{i}$;

ii) $\hat\rho_{\varphi}=\hat{w}_{1} \dots \hat {w} _ {q} $.

We choose numbers $n_i\in\mathbb N$ so that $\bar\varphi^{-n_q}(B_{j(q)})\subset\dots\subset\bar\varphi^{-n_1}(B_{j(1)})$. Let $\bar K_i=\bar\varphi^{-n_{j(i)}}(B_{j(i)})\setminus \bar\varphi^{-n_{j(i)}-1}(B_{j(i)})$. 
Denote by $\bar w_{i,t}:\mathbb R^2\to\mathbb R^2$ diffeomorphism coinciding with $(p_\varphi|_{\bar K_{i}})^{-1}\hat w_{i,t}p_\varphi$ on $\bar K_i$ and identical outside  $\bar K_i$. Let $\bar w_t=\bar w_{1,t}\dots,\bar w_{q,t}$ and $\bar\nu_{\varphi,t}=\bar\varphi\bar w^{-1}_t.$ Denote by $\nu_{\varphi,t}:M^2\to M^2$ a diffeomorphism coinciding with $\psi_{\alpha_1}^{-1}\bar \nu_{\varphi,t} \psi_{\alpha_1}$ on $B_\varphi$ and coinciding with $\varphi$ outside $B_\varphi$. Then  $\nu_{\varphi,1}\in\Phi_{\varphi_1}$ since $\bar\nu_{\varphi,1}$ coincides with $\bar\varphi_1$ on some disk $B_{\nu_{\varphi,1}}=\bar\varphi^{-m}(B_\varphi)\subset\bar\varphi^{-n_{j(q)}-1}(B_{j(q)})$. By construction $\bar\rho_{\nu_{\varphi,1}}=\bar\varphi_1^{-m}\bar\rho_{\varphi}(\bar\varphi\bar w^{-1}_1)^{m}$ and therefore $\hat\rho_{\nu_{\varphi,1}}(\hat x)=\hat x,\,\hat x\in \mathbb T^2$. 

Thus, the diffeomorphism $\rho_{\nu_{\varphi,1}}$ in a neighborhood of the point $\alpha_1$ coincides with some degree of the diffeomorphism $\varphi$ and, therefore, smoothly extends to the point $\alpha_1$. 

In subcase b), construct an arc without bifurcations $\mu_{\varphi,t}:M^2\to M^2$ such that $\mu_{\varphi,0}=\varphi$, $\mu_{\varphi,1}\in \Phi_{\varphi_1}$ and  $n_{\mu_{\varphi,1}}=0$. Then the bifurcation free arc $$H_{\varphi,\alpha_1,t}=\mu_{\varphi,t}*\nu_{\mu_{\varphi,1},t}*(\varrho_{\nu_{\mu_{\varphi,1},1},t}^{-1}\nu_{{\mu_{\varphi,1},1}}\varrho_{\nu_{\mu_{\varphi,1},1},t})$$ connects the diffeomorphism $\varphi$ with the diffeomorphism $\varphi_1$. 

We describe the construction of the arc $\mu_{\varphi,t}$.
 
Introduce a coordinates $r,\phi$ on the disk $B_\varphi$ in which the curve $\bar\varphi^{-k}(\partial B_\varphi),\,k\in\mathbb N$ has the form $r=2^{-k},\,\phi\in[0,2\pi)$. Define the diffeomorphism $\bar\theta_t:\mathbb R^2\to \mathbb R^2$ so that $\bar\theta_t(O)=O$ and $\bar{\theta}_t(r,\phi)=r e^{i(\phi+\phi_t(r))}$, where $\phi_t(r)$ is a smooth monotone function equal to $-2n_{\varphi}\pi t$ for $r\leqslant\frac{1}{2}$ and equal to $0$ for $r\geqslant{1}$. Let $ \bar\mu_{\varphi,t}=\bar\varphi\bar\theta_t$ and denote by $\mu_{\varphi,t}:M^2\to M^2$ diffeomorphism coinciding with $\psi_{\alpha_1}^{-1}\bar \mu_{\varphi,t} \psi_{\alpha_1}$ on $B_\varphi$ and coinciding with $\varphi$ outside $B_\varphi$. Then $\mu_{\varphi,1}\in\Phi_{\varphi_1}$ since $\bar\mu_{\varphi,1}$ is the same as $\bar\varphi_1$ on disk $B_{\mu_{\varphi,1}}=\bar\varphi^{-1}(B_\varphi)$. By construction $\bar\rho_{\mu_{\varphi,1}}=\bar\varphi_1^{-1}\bar\rho_{\varphi}\bar\varphi\bar\theta_1$ and therefore $n_{\mu_{\varphi,1}}=0$.
\end{demo}

\section{Construction an arc $H_{J,t}$}

\subsection{Construction of auxiliary functions}

In this section, we construct model functions that will later be used to construct a stable arc. The construction is based on the principle of gluing infinitely smooth functions by means of the following {\it sigmoid function}.

\begin{figure}[h]
\centerline{\includegraphics[width=8 true cm]{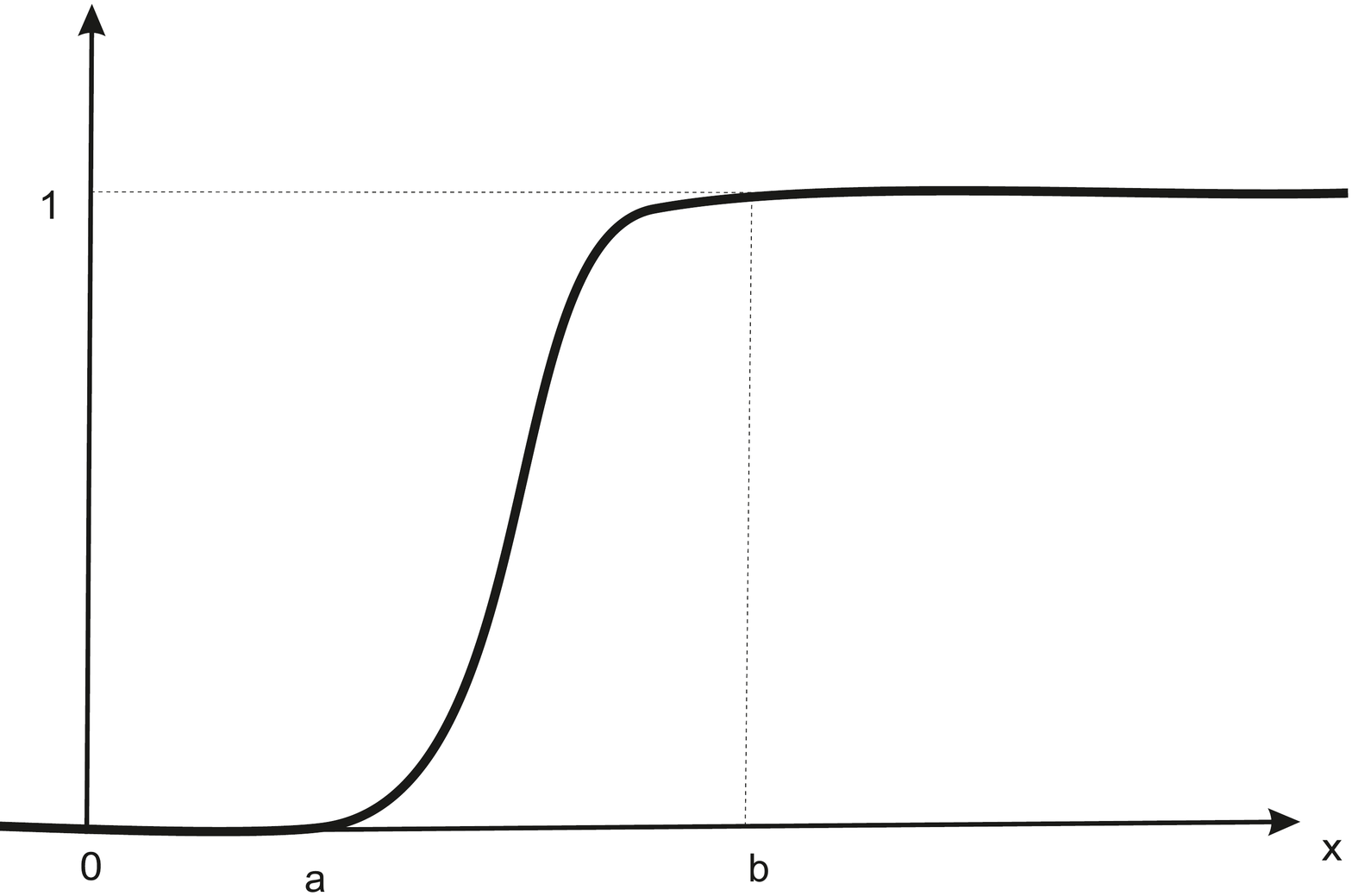}}
\caption{Graph of the sigmoid function}\label{si}
\end{figure}

Let $a<b$ and $\delta_{_{a;b}}:\mathbb R\to[0,1]$ sigmoid function, defined by the formula (Fig. \ref{si}) $$\delta_{_{a;b}}(x) =\begin{cases} 0,&x\leqslant a,\cr\frac{1}{1+\exp\left(\frac{(a+b)/2-x}{{(x-a)}^2{(x-b)}^2}\right)},&a<x<b,\cr 1,&x\geqslant b.\end{cases}$$ 

\begin{figure}[h!]
\centerline{\includegraphics[width=8 cm]{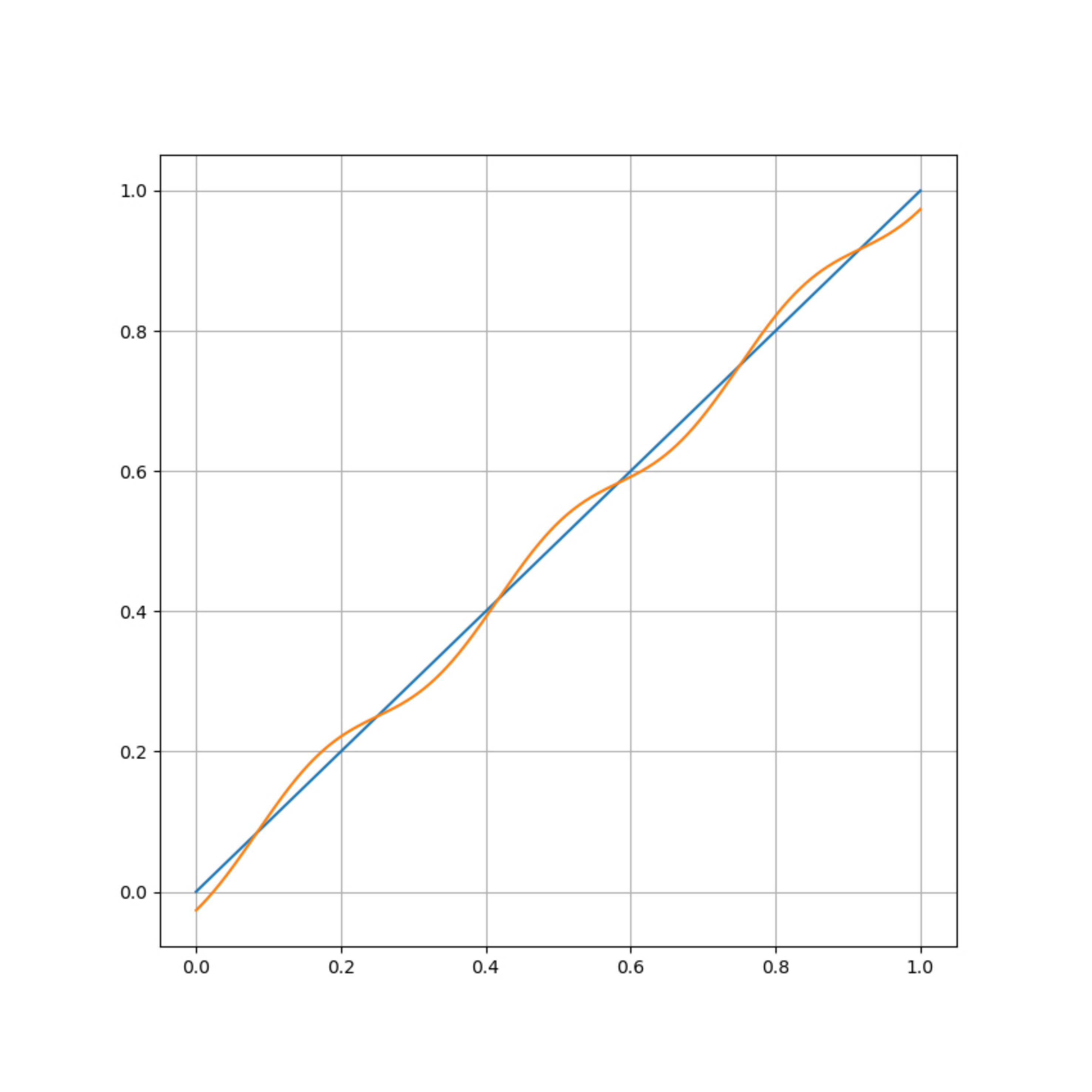}}
\caption{\small Graph of the map $\bar \phi_1 (x)$}\label{phi_1}
\end{figure}

Define the function $\bar \phi_1:\mathbb R\to\mathbb R$ by the formula (Fig. \ref{phi_1}).
$$\bar \phi_1(x)= x-\frac{1}{12\pi}sin\left(6\pi\left (x-\frac{1}{4}\right)\right).$$

\begin{figure}[h]
\centerline{\includegraphics[width=10 cm]{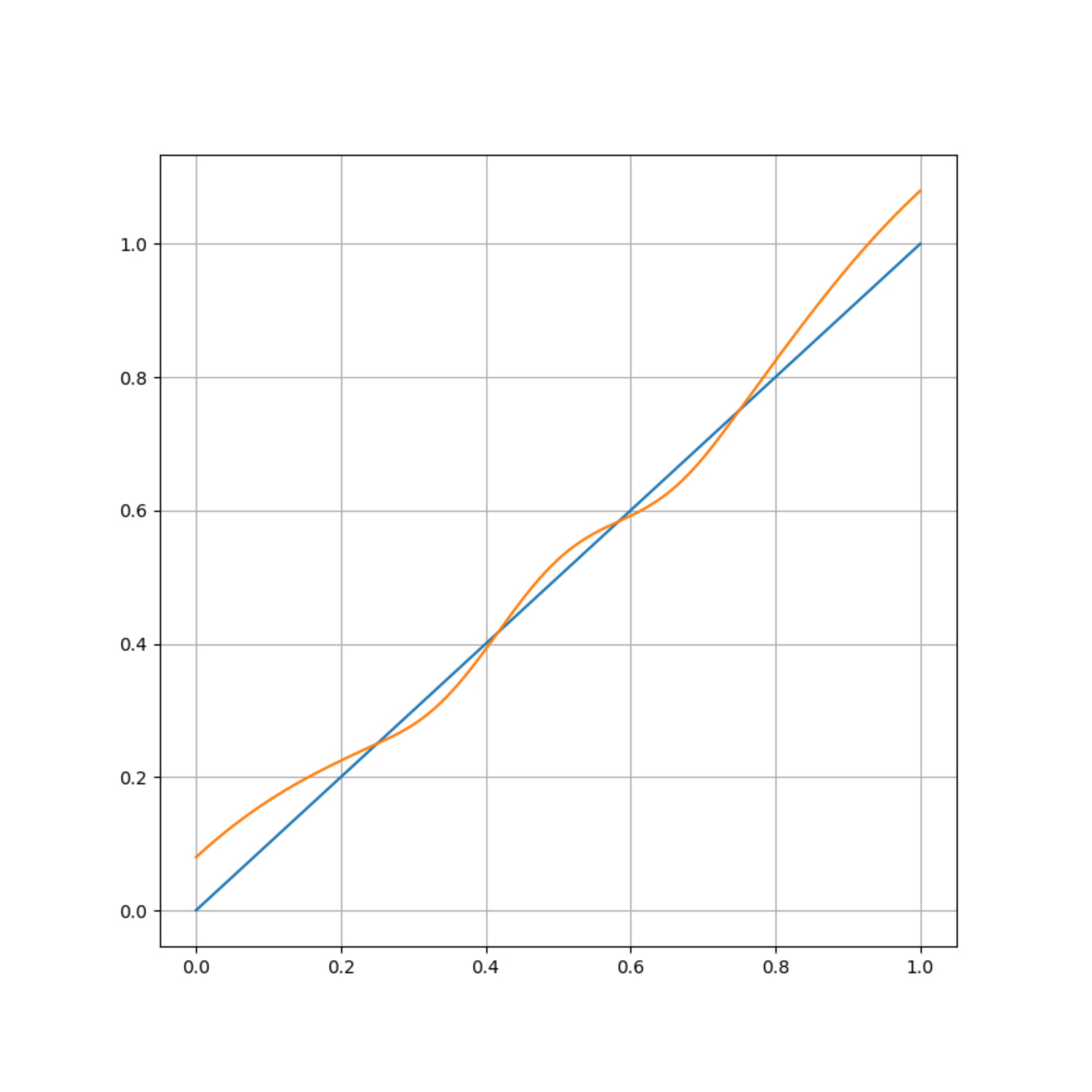}}
\caption{\small Graph of the map $\bar g_1(x)$}\label{g1}
\end{figure}

Define the function $\bar g_1:\mathbb R\to\mathbb R$ by the formula (Fig. \ref{g1})
$$\bar g_1 (x) = 
 \begin{cases}
 \bar {\phi}_0(x),\, 0\leqslant x \leqslant 0,26,\\
 (1-\delta_{_{0,26; 0,27}}(x))\bar \phi_0(x)+\delta_{_{0,26; 0,27}}(x) \bar \phi_1(x),\, 0,26< x <0,27,\\
 \bar {\phi}_1(x),\, 0,27 \leqslant x \leqslant 0,76,\\
 (1-\delta_{_{0,76; 0,77}}(x))\bar \phi_1(x)+\delta_{_{0,76; 0,77}}(x) \bar \phi_0(x),\, 0,76< x <0,77,\\
  \bar {\phi}_0(x),\, 0,77 \leqslant x \leqslant 1\\
 \end{cases}.$$ 

\begin{figure}[h!]
\centerline{\includegraphics[width=8 cm]{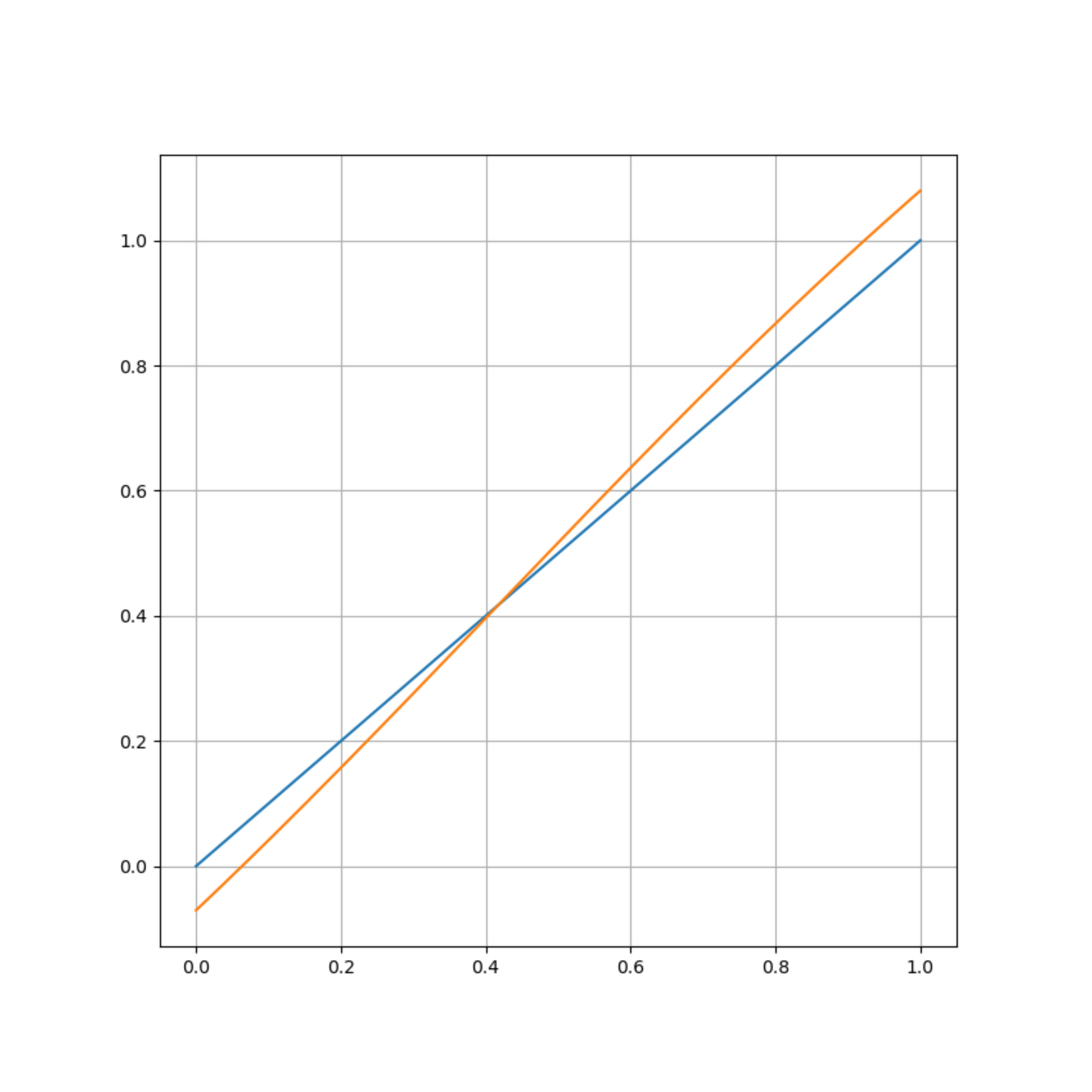}}
\caption{\small Graph of the map $\bar \phi_2 (x)$}\label{phi_2}
\end{figure}

Define the function $\bar \phi_2:\mathbb R\to\mathbb R$ by the formula (Fig. \ref{phi_2}).
$$\bar \phi_2(x)= x+\frac{1}{4\pi}sin\left(\frac{5}{6}\pi\left (x-\frac{5}{12}\right)\right).$$

\begin{figure}[h!]
\centerline{\includegraphics[width=8 cm]{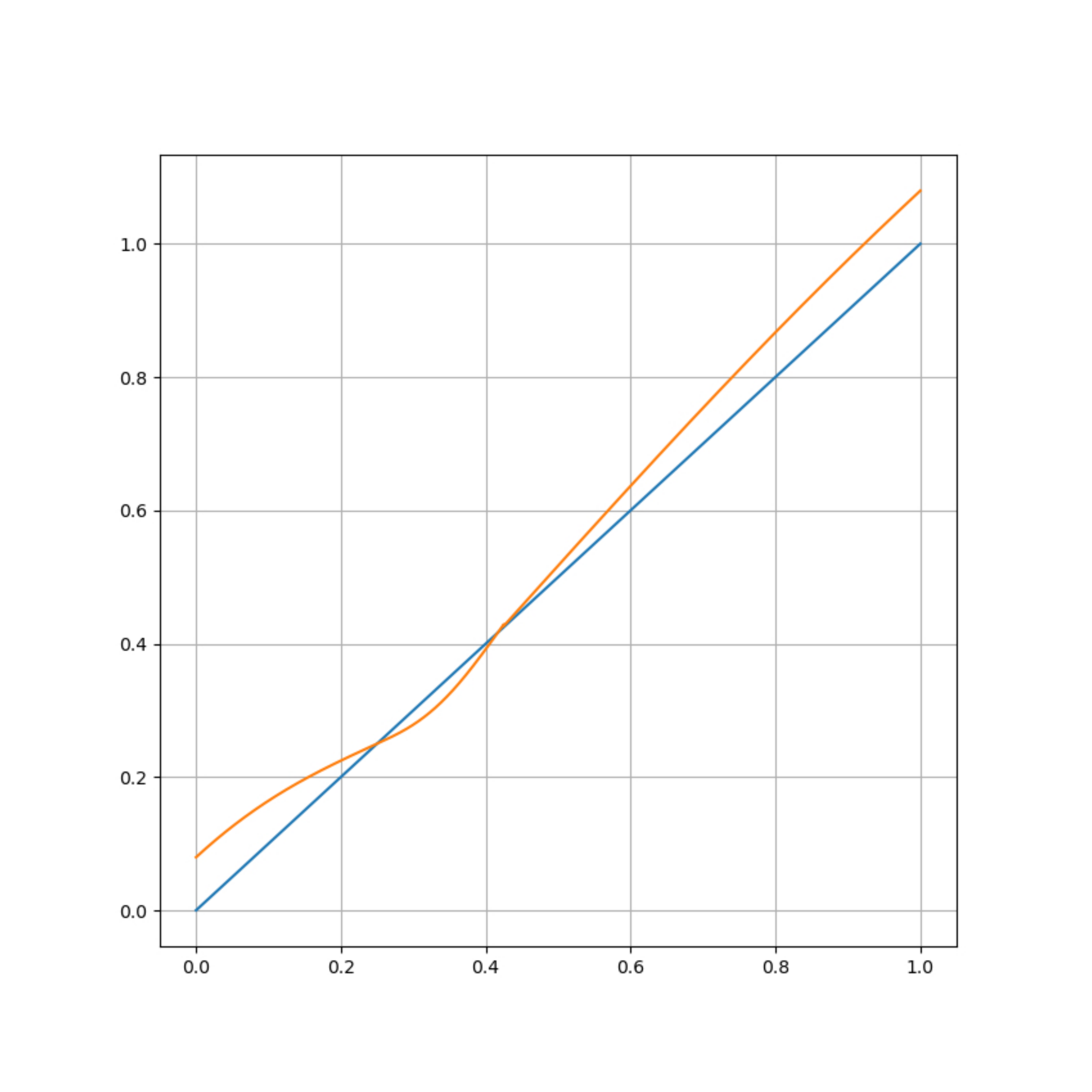}}
\caption{\small Graph of the map $\bar g_2 (x)$}\label{g_2}
\end{figure}

Define the function $\bar g_2:\mathbb R\to\mathbb R$ by the formula (Fig. \ref{g_2}) 
$$\bar g_2 (x) = 
 \begin{cases}
 \bar {g}_1(x),\, 0\leqslant x \leqslant 0,42,\\
(1-\delta_{_{0,42; 0,43}}(x))\bar g_1(x)+ \delta_{_{0,42; 0,43}}(x) \bar \phi_2(x),\, 0,42< x <0,43,\\
  \bar {\phi}_2(x),\, 0,43 \leqslant x \leqslant 0,98,\\
(1-\delta_{_{0,98; 0,99}}(x))\bar \phi_2(x)+ \delta_{_{0,98; 0,99}}(x) \bar g_1(x),\, 0,98< x <0,99,\\
  \bar {g}_1(x),\, 0,99 \leqslant x \leqslant 1\\
 \end{cases}.$$ 

\subsection{Construction the model arcs}

In this section, we will construct arcs, which are the main components that make up an arc $H_{J,t}$. 

For $n\in\mathbb Z$ let $J_n=\begin{pmatrix}
  1& 0\\
  n& 1
\end{pmatrix}.$

\begin{lemma}\label{11} 
The diffeomorphism $f_0$ is connected with the diffeomorphism $f_{J_1}$ by a stable arc $H_{0,1,t}$ with  typically passing non-critical two saddle-node bifurcations.
\end{lemma}
\begin{demo} In this proof, the bar-free mappings are projections on $\mathbb S^1$ by $ \pi $ of the bar-mappings given on the line $\mathbb R$. The stable arc $H_{0,1,t}$, connecting the diffeomorphism $f_0$ with the diffeomorphism $f_{J_1}$ product of the arcs $\Gamma^1_t,\,\Gamma^2_t$, constructed in step 1 and step 2 below and the arc $H_{\Gamma^2_1,t}$.

{\bf Step 1. First saddle-node bifurcation.}

{\bf 1. The birth of a saddle-node point.}

We start with the diffeomorphism $f_0: \mathbb T^2\to\mathbb T^2$, defined by the formula:
$$f_0(z,w)=(\phi_0(z),\phi_0(w)),\,z,w\in\mathbb S^1.$$
Let $$\bar\eta^1_t(x)=(1-t)\bar{\phi}_0(x)+t\bar g_1(x),\,x\in\mathbb R,\,t\in[0,1]$$ and 
$$\bar\eta^1_{t,\tau}(x)=(1-\tau)\bar{\eta}^1_t(x)+\tau\bar \phi_0(x),\,x\in\mathbb R,\,t\in[0,1],\,\tau\in[0,1].$$  Define a smooth arc $H^1_t: \mathbb T^2\to\mathbb T^2$, $t\in [0,1]$ by the formula:
$$H^1_t(w,z)=\begin{cases}
(\phi_0(\pi(x)),\eta^1_{t,|8x-2|}(z)),\,x\in\left(\frac18,\frac38\right),\,z\in\mathbb S^1,\\
f_0(\pi(x),z),\,x\in\left(-\frac{5}{8},\frac18\right),\,z\in\mathbb S^1\\
 \end{cases}.$$ 

For $t=\frac{3}{4}$, the diffeomorphism $H^1_\frac34$ has a saddle-node point $p=(N,\pi(0))$, whose stable manifold is diffeomorphic to a half-plane whose boundary is arc $\gamma_p$ (Fig. \ref{Ht}). 
\begin{figure}[h]
\begin{minipage}[h]{0.23\linewidth}
\center{\includegraphics[width=1\linewidth]{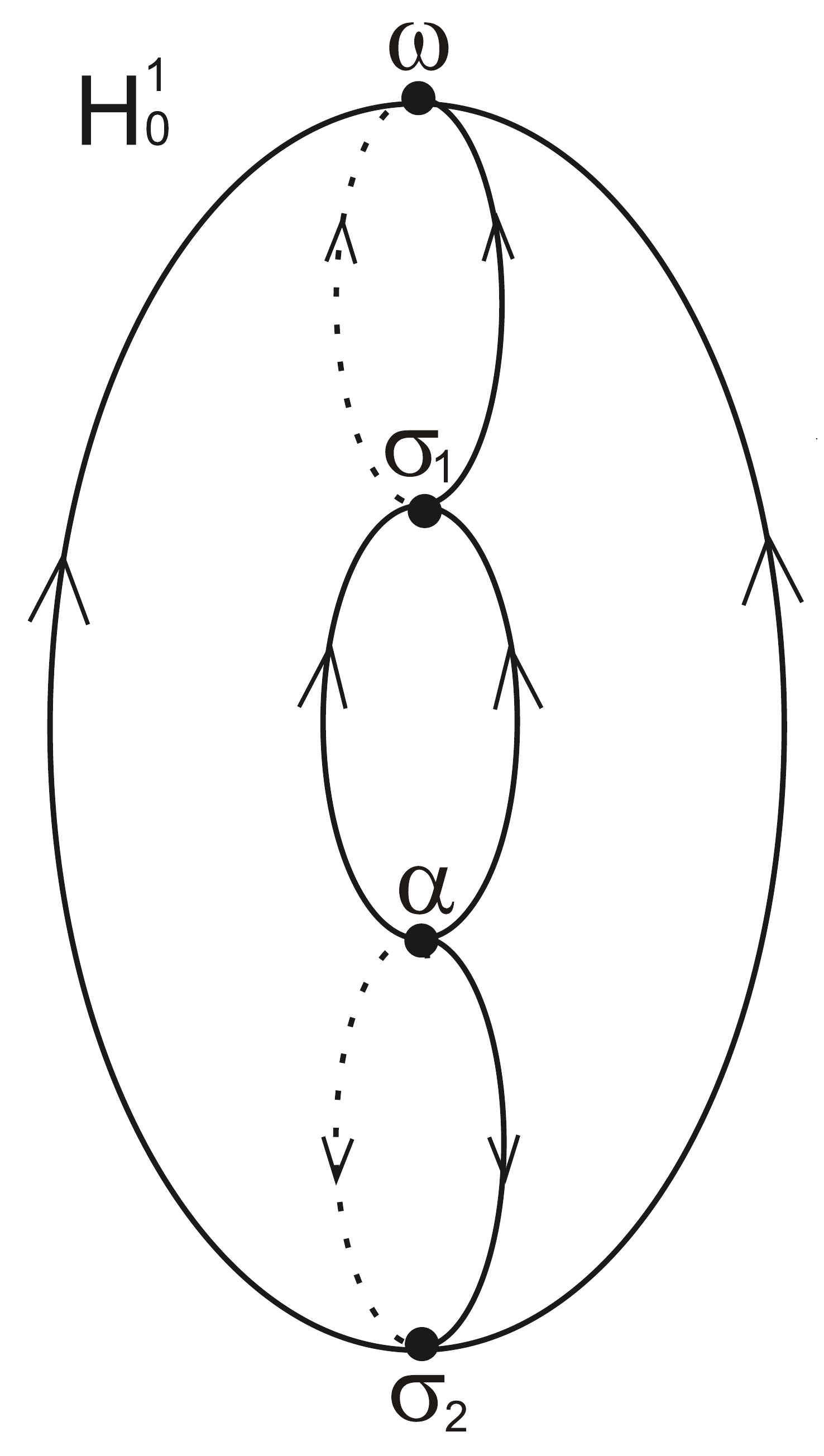}}
\end{minipage}
\hfill
\begin{minipage}[h]{0.23\linewidth}
\center{\includegraphics[width=1\linewidth]{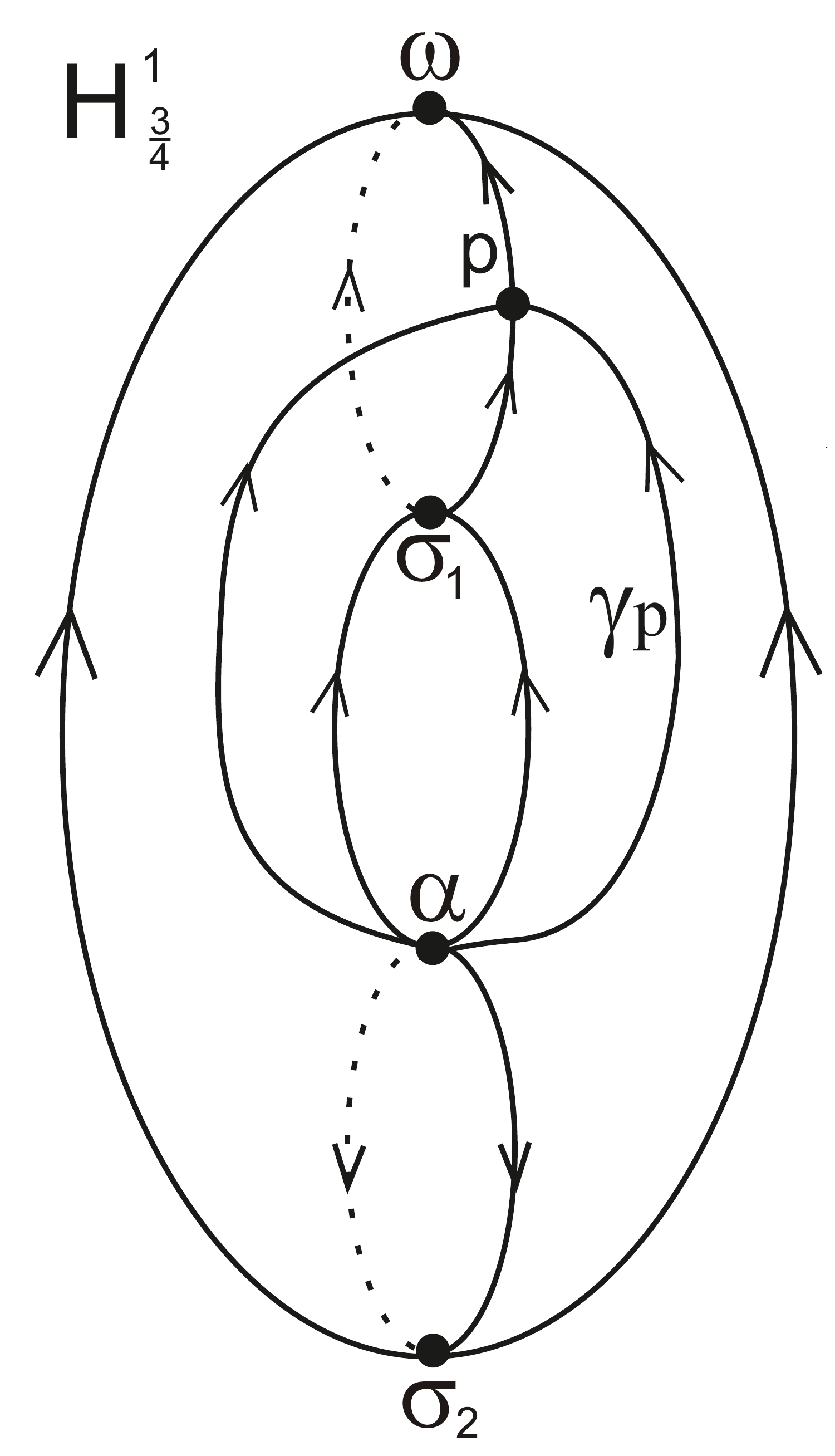}}
\end{minipage}
\hfill 
\begin{minipage}[h]{0.23\linewidth}
\center{\includegraphics[width=1\linewidth]{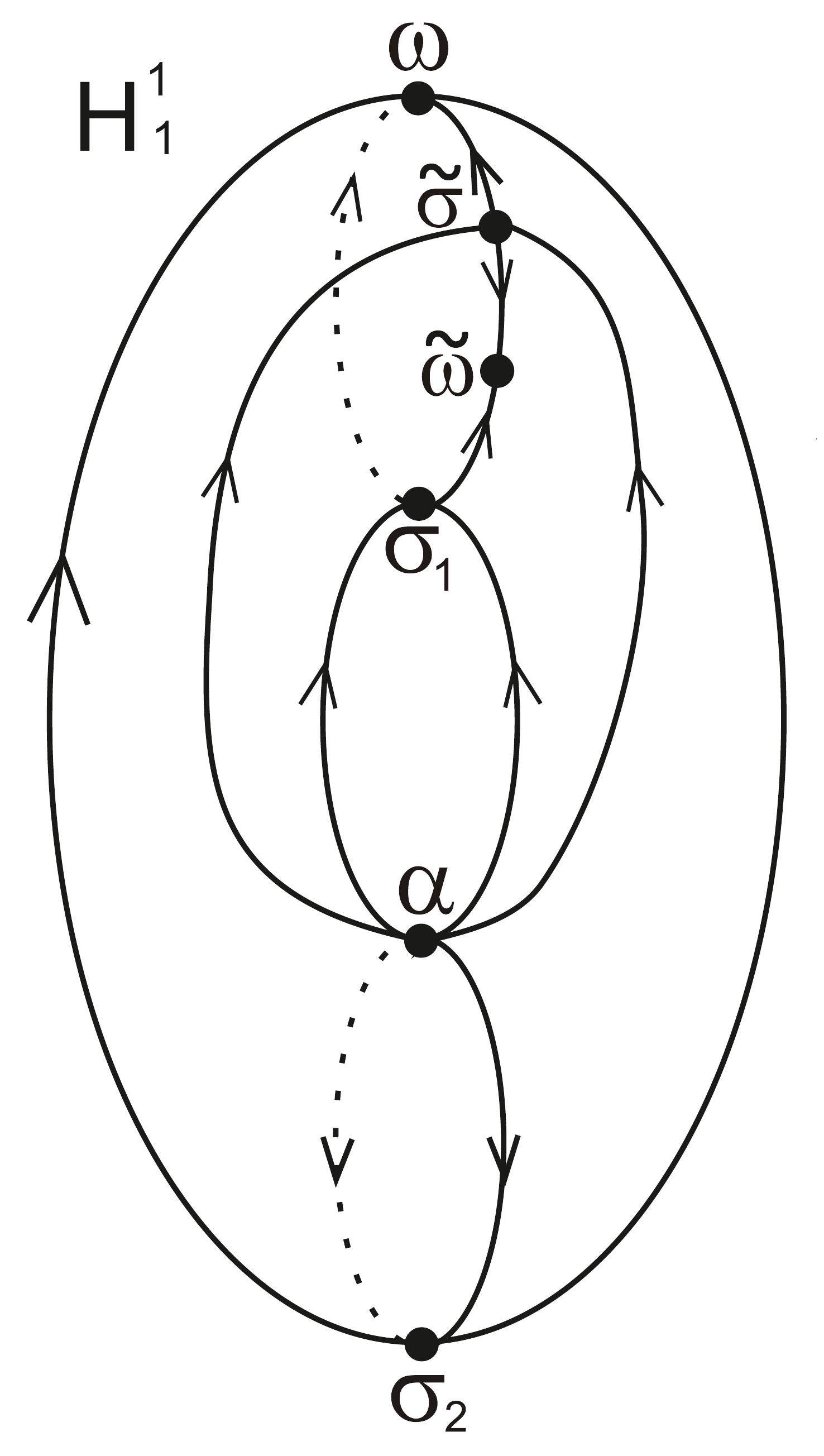}}
\end{minipage}
\hfill 
\caption{Isotopy $H^1_t$ on the torus}
\label{Ht}
\end{figure}

{\bf 2. Rotation of the separatrix of the saddle $\sigma_2$.}

Consider the fundamental domain $K=\left[\pi(0),\pi\left(\frac{1}{4\pi}\right)\right] \times \mathbb S^1$ restriction of the diffeomorphism $f_0$ to $V=\left[\pi\left(-\frac{1}{4}\right),\pi\left(\frac{1}{4}\right)\right] \times \mathbb S^1 $. 
Let $\hat V=V/f_0$. Then $\hat V$ is a torus got from $K$ by identifying the boundaries with the map $f_0$. Denote by $q:V\to\hat V$ the natural projection. Let  $\hat \gamma_{2}=q(W^u_{\sigma_2}\cap V)$ and $\hat \gamma_{1}=q(W^s_{\sigma_1}\cap V)$. Since the diffeomorphism $H_t$ for all $t\in[0,1]$ coincides with $f_0$ on the annulus $\left[\pi(-\frac{1}{4}),\pi\left(\frac{1}{8}\right)\right]\times \mathbb S^1 $, then the circle $\hat\gamma_p=q(\gamma_p\cap K)$ is defined correctly.

Let $W=\left[\pi\left(-\frac{1}{4}\right);\pi\left(\frac{1}{4}\right)\right] \times \left[\pi\left(-\frac{1}{4}\right);\pi\left(\frac{1}{4}\right)\right]$ and $\hat W=p(W)$. By construction, the circle $\hat\gamma_p$ divides the annulus $\hat W$ into two annulus, the closures of which are denoted by$\hat W_1,\hat W_2$, assuming that $\hat\gamma_{1}\subset\hat W_1$ and $\hat\gamma_{2}\subset\hat W_2$ (Fig. \ref{zeta}).

Choose a circle $\hat\gamma\subset int\,\hat W_1$ that is not homotopic to zero. According to \cite{Rolf}, there exists a diffeomorphism  $\hat h_1:\hat V\to\hat V$ smoothly isotopic to the identity such that $\hat h_1 (\hat\gamma_2)=\hat\gamma$. 
 
For $x_i\in [-\frac{1}{4};0]$ let $K_{i}=[ \pi(x_i);(\pi(\bar\phi_0^{-1}(x_i))] \times \mathbb S^1$. Choose an open cover $D=\{D_1,\dots, D_{k_1}\}$ of the torus $\mathbb{T}^2$ such that the connected component $\bar D_i$ of the set $q^{-1}(D_i)$ is a subset of $K_{i} $ for some $x_i<\bar\phi_0^{-1}(x_{i-1})$. According to the fragmentation lemma \cite{Ba}  there exist smoothly isotopic to the identity diffeomorphisms $\hat{w}_ {1}, \dots, \hat {w} _ {k_1}: \mathbb {T}^2\to\mathbb{T}^2$ with the following properties:

i) for each $i\in\{1,\dots, k_1\}$  there exists smooth isotopy $\{\hat{w}_{i,t}\}$ which is the identity outside $D_{i}$ and which joins the identity and $\hat{w}_{i}$

ii) $\hat{h}_1=\hat{w}_{1} \dots \hat {w} _ {k_1}$.

Let ${w}_{i,t}:\mathbb{R}^2\to\mathbb{R}^2$ diffeomorphism that coincides with 
$(q\vert_{K_{i}})^{-1}\hat{w}_{i,t}q$ на $K_{i}$ and coincides with the identity map outside $K_{i}$. 
Let $${\zeta}_{t}={w}_{1,t}\dots {w}_{k_1,t}f_0,\,\,\,\,G^1_t=\begin{cases}
\zeta_{2t},\,0\leqslant t<\frac12,\\
\zeta_1,\,\frac12\leqslant t\leqslant 1\\
 \end{cases}.$$ 
\begin{figure}[h]
\center{\includegraphics[width=.5\linewidth]{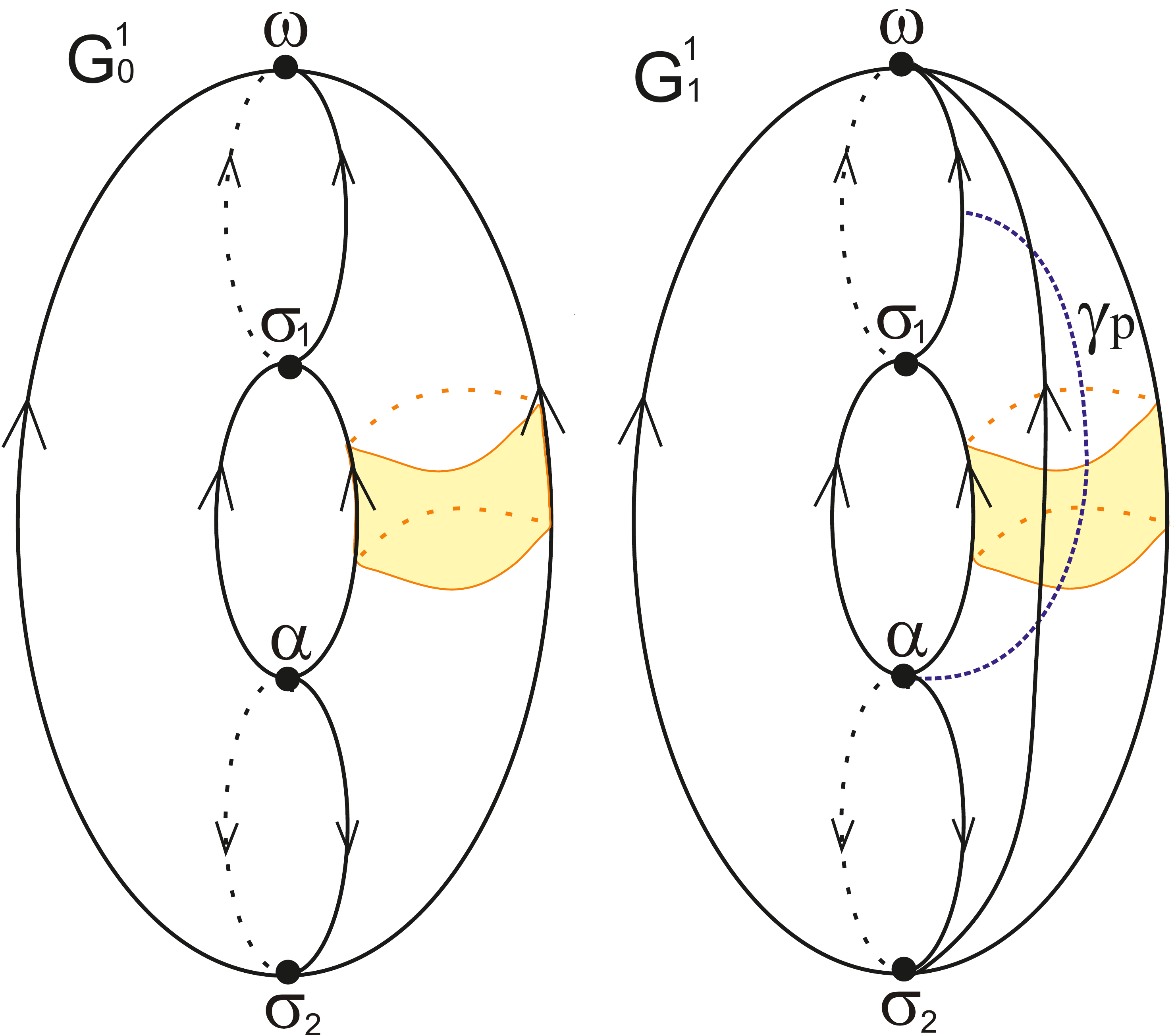}}
\caption{Isotopy ${G}^1_{t}$ on the torus}
\label{zeta}
\end{figure}

{\bf 3. Сombining isotopies $H^1_t$ and $G^1_t$.}

Define a smooth arc $\Gamma^1_t :\mathbb T^2\to\mathbb T^2$, $t\in [0,1]$ by the formula (Fig. \ref{Gam}):
$$\Gamma^1_t(z,w)=\begin{cases}
H^1_t(\pi(x),w),\,x\in\left(\frac18,\frac38\right),\,z\in\mathbb S^1,\\
G^1_t(\pi(x),w),\,x\in\left(-\frac14,0\right),\,z\in\mathbb S^1,\\
f_0(\pi(x),w),\,x\in\left[-\frac{5}{8},-\frac14\right]\cup [0;\frac38],\,z\in\mathbb S^1\\
 \end{cases}.$$ 

\begin{figure}[h]
\begin{minipage}[h]{0.23\linewidth}
\center{\includegraphics[width=1\linewidth]{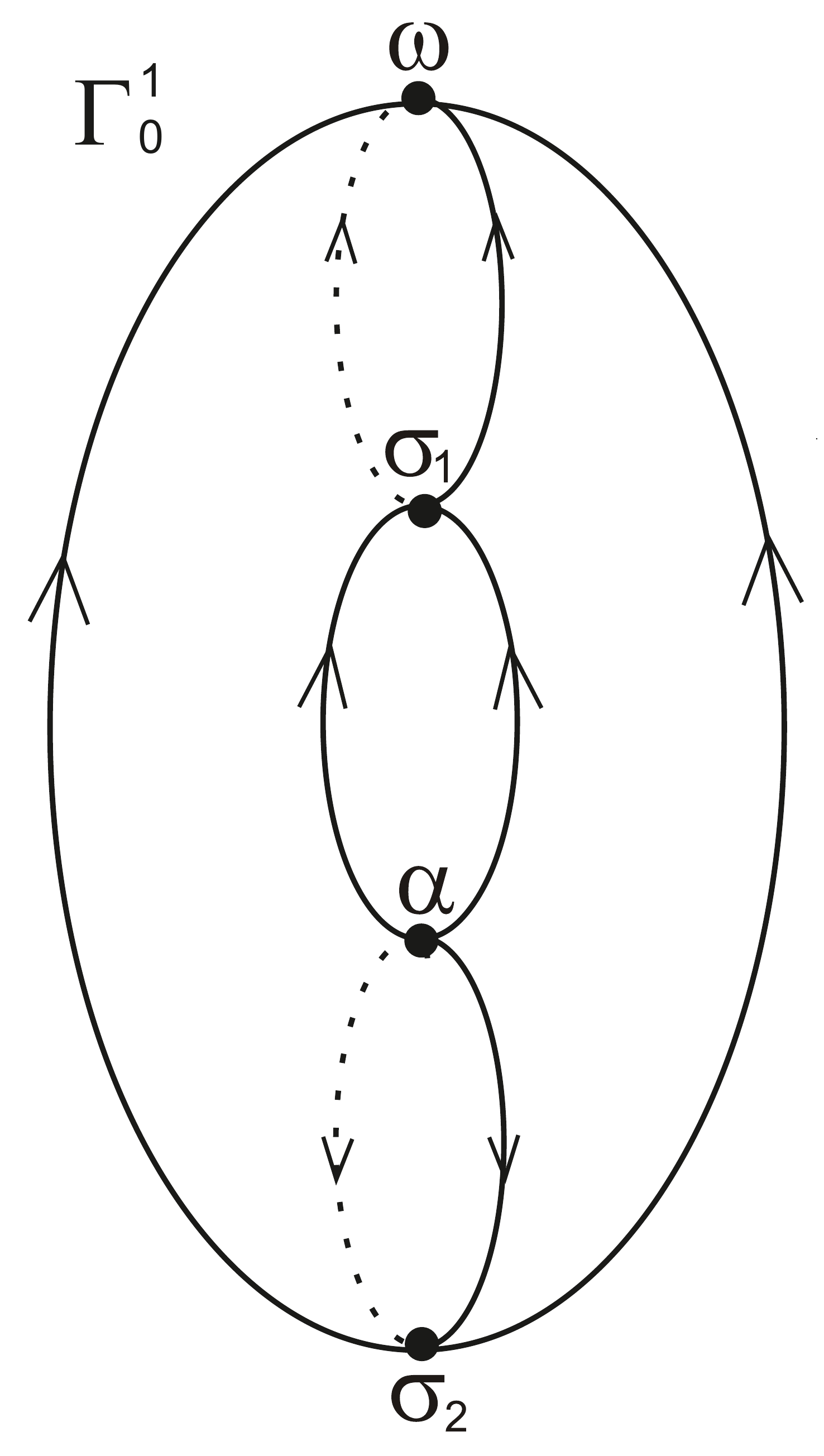}}
\end{minipage}
\hfill
\begin{minipage}[h]{0.23\linewidth}
\center{\includegraphics[width=1\linewidth]{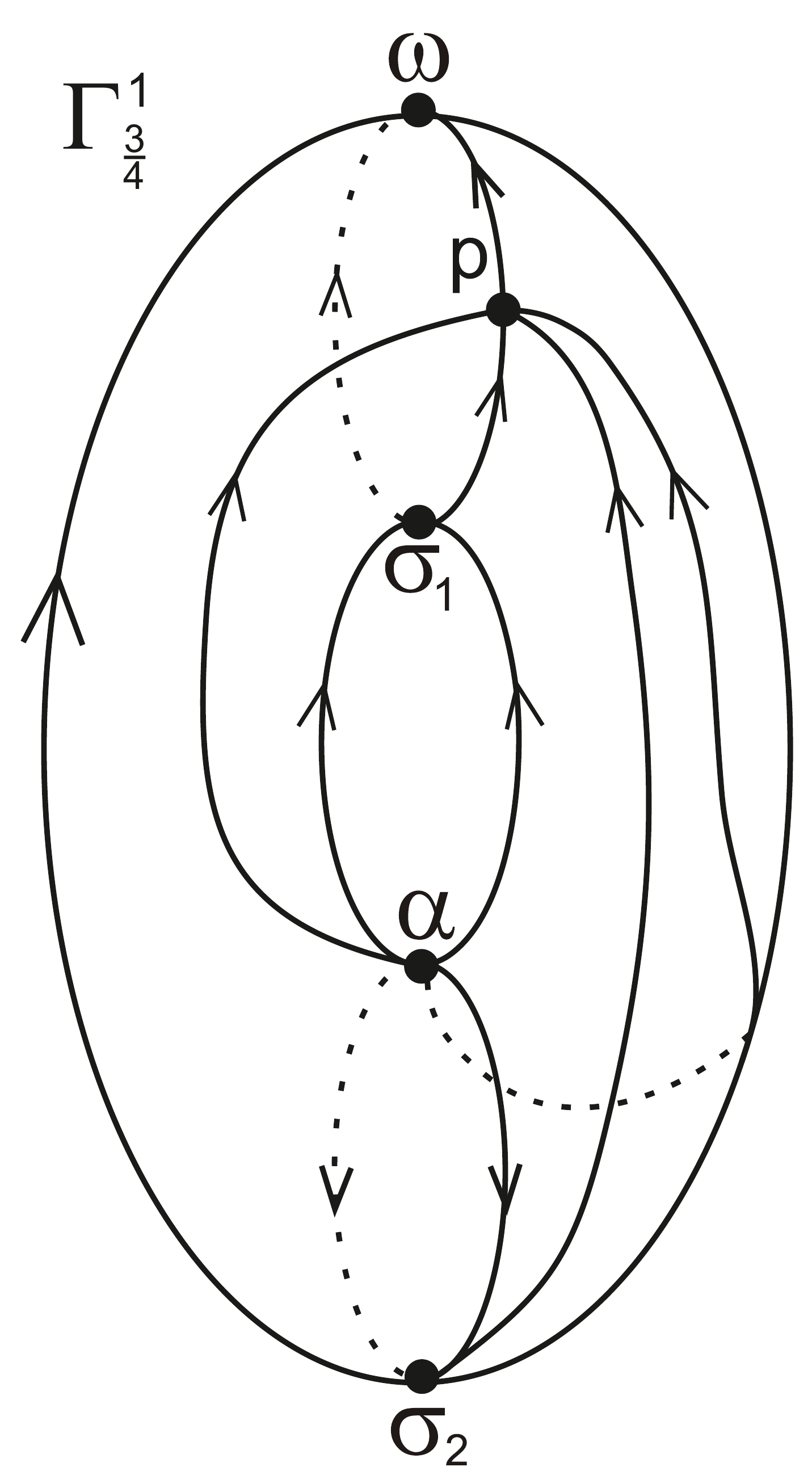}}
\end{minipage}
\hfill 
\begin{minipage}[h]{0.23\linewidth}
\center{\includegraphics[width=1\linewidth]{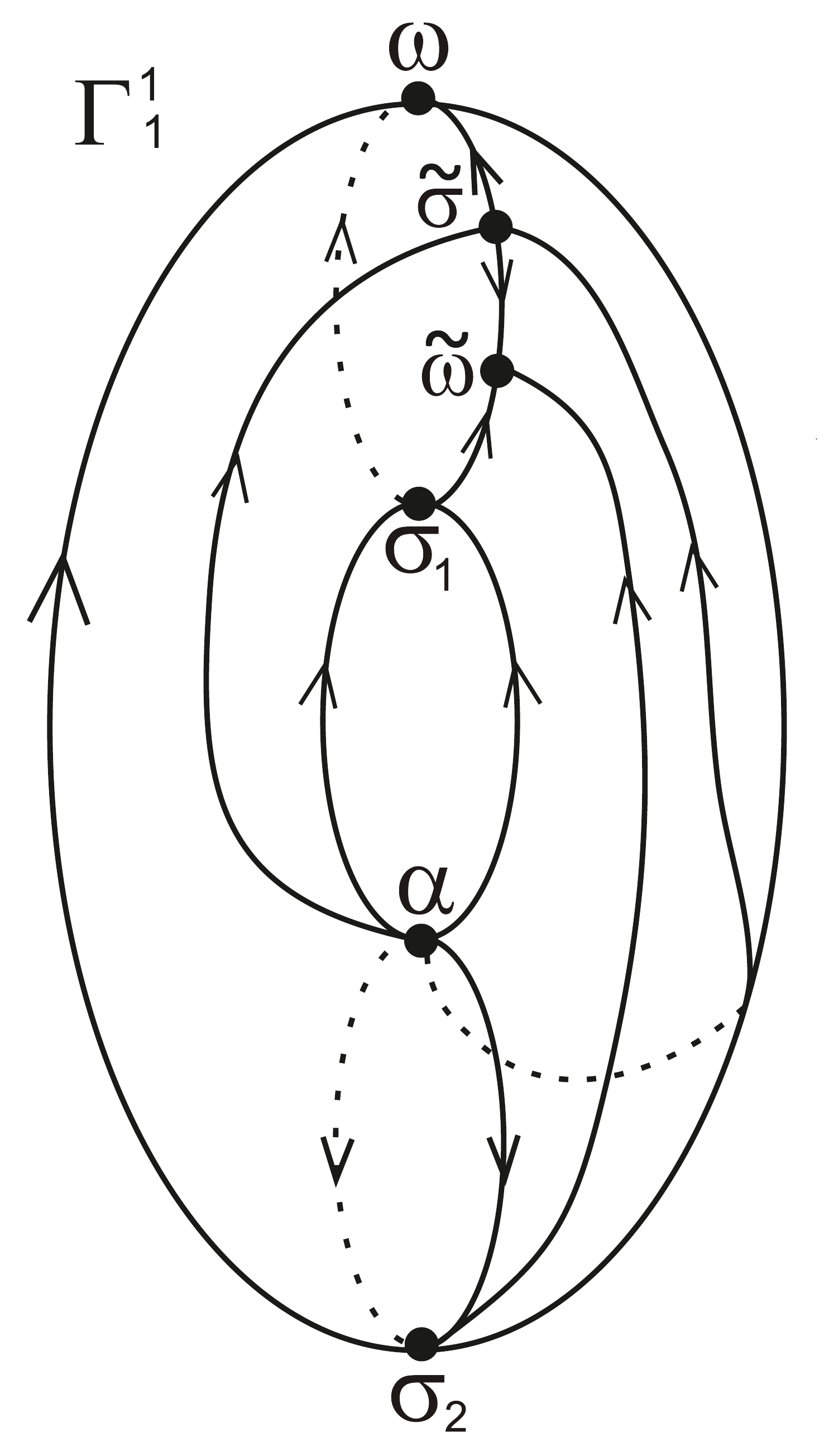}}
\end{minipage}
\hfill 
\caption{Isotopy $\Gamma^1_t$ on the torus}
\label{Gam}
\end{figure}

{\bf Step 2. Second saddle-node bifurcation.}

{\bf 1. Merging saddle and node points}

For all $t\in [0;1]$ let $\bar{\eta}^2_t(x)=t\bar{g_2}(x)+(1-t)\bar{g_1}(x), x\in\mathbb{R}$ and $$\bar\eta^2_{t,\tau}(x)=(1-\tau)\bar{\eta}^2_t(x)+\tau\bar \phi_0(x),\,x\in\mathbb R,\,t\in[0,1],\,\tau\in[0,1].$$ 

Define a smooth arc $H^2_t: \mathbb T^2\to\mathbb T^2$, $t\in [0,1]$ by the formula:
$$H^2_t(w,z)=\begin{cases}
(\phi_0(\pi(x)),\eta^2_{t,|8x-2|}(z)),\,x\in\left(\frac18,\frac38\right),\,z\in\mathbb S^1,\\
\Gamma_1(\pi(x),z),\,x\in\left(-\frac{5}{8},\frac18\right),\,z\in\mathbb S^1\\
 \end{cases}.$$ 

The arc $H^2_t$ realizes the merging of the sink $\tilde\omega$ and the saddle $\sigma_1$ into the saddle-node point $\tilde p$ and its further disappearance. Denote by$\beta_{\tilde p}$ the boundary of the stable manifold of a saddle-node $\tilde p$.

{\bf 2. Rotation of the separatrix of the saddle $\sigma_2$.}

Since the diffeomorphism $H^2_t$ for all $t\in[0,1]$ coincides with $f_0$ on the annulus $K$, then the circles $\hat \beta_{2}=q(W^u_{\sigma_2}\cap K)$, $\hat \beta_{1}=q(W^s_{\tilde {\sigma}}\cap K)$ and $\hat\beta_{\tilde {p}}=q(\beta_{\tilde {p}}\cap K) is defined correctly$.

 $\hat\gamma_p=q(\gamma_p\cap K)$ .

Let $\hat W_3=\hat V\setminus \hat W$. Choose a circle $\hat\beta\subset int\,\hat W_3$ that is not homotopic to zero. According to \cite{Rolf}, there exists a diffeomorphism $\hat h_2:\hat V\to\hat V$ smoothly isotopic to the identity such that $\hat h_2 (\hat\beta_{2})=\hat\beta$ and $\hat h_2 (\hat\beta_{1})=\hat\beta_1$. 

Choose an open cover $U=\{U_1,\dots, U_{k_2}\}$ of the torus $\mathbb{T}^2$ such that the connected component $\bar U_i$ of the set $q^{-1}(U_i)$ is a subset of $K_{i} $ for some $x_i<\bar\phi_0^{-1}(x_{i-1})$. According to the lemma of fragmentation \cite{Ba} there exist smoothly isotopic to the identity diffeomorphisms $\hat{v}_ {1}, \dots, \hat {v}_{k_2}: \mathbb {T}^2\to\mathbb{T}^2$ with the following properties:

i) for each $i\in\{1,\dots, k_2\}$ there exist a smooth isotopy $\{\hat{v}_{i,t}\}$ which is the identity outside $U_{i}$ and which joins the identity and $\hat{v}_{i}$;

ii) $\hat{h}_2=\hat{v}_{1} \dots \hat {v} _ {k_2} $.

Let ${v}_{i,t}:\mathbb{R}^2\to\mathbb{R}^2$ diffeomorphism that coincides with $(q\vert_{K_{i}})^{-1}\hat{v}_{i,t}q$ on $K_{i}$ and coincides with the identity map outside $K_{i}$. 
Let $${\xi}_{t}={v}_{1,t}\dots {v}_{k_2,t}\Gamma_1,\,\,\,\,G^2_t=\begin{cases}
\xi_{2t},\,0\leqslant t<\frac12,\\
\xi_1,\,\frac12\leqslant t\leqslant 1\\
 \end{cases}.$$ 

{\bf 3. Сombining isotopies $H^2_t$ and $G^2_t$.}

Define a smooth arc $\Gamma^2_t :\mathbb T^2\to\mathbb T^2$, $t\in [0,1]$ by the formula (Fig. \ref{Gam}):
$$\Gamma^2_t(z,w)=\begin{cases}
H^2_t(\pi(x),w),\,x\in\left(\frac18,\frac38\right),\,z\in\mathbb S^1,\\
G^2_t(\pi(x),w),\,x\in\left(-\frac14,0\right),\,z\in\mathbb S^1,\\
f_0(\pi(x),w),\,x\in\left[-\frac{5}{8},-\frac14\right]\cup [0;\frac38],\,z\in\mathbb S^1\\
 \end{cases}.$$ 

\begin{figure}[h]
\begin{minipage}[h]{0.32\linewidth}
\center{\includegraphics[width=0.9\linewidth]{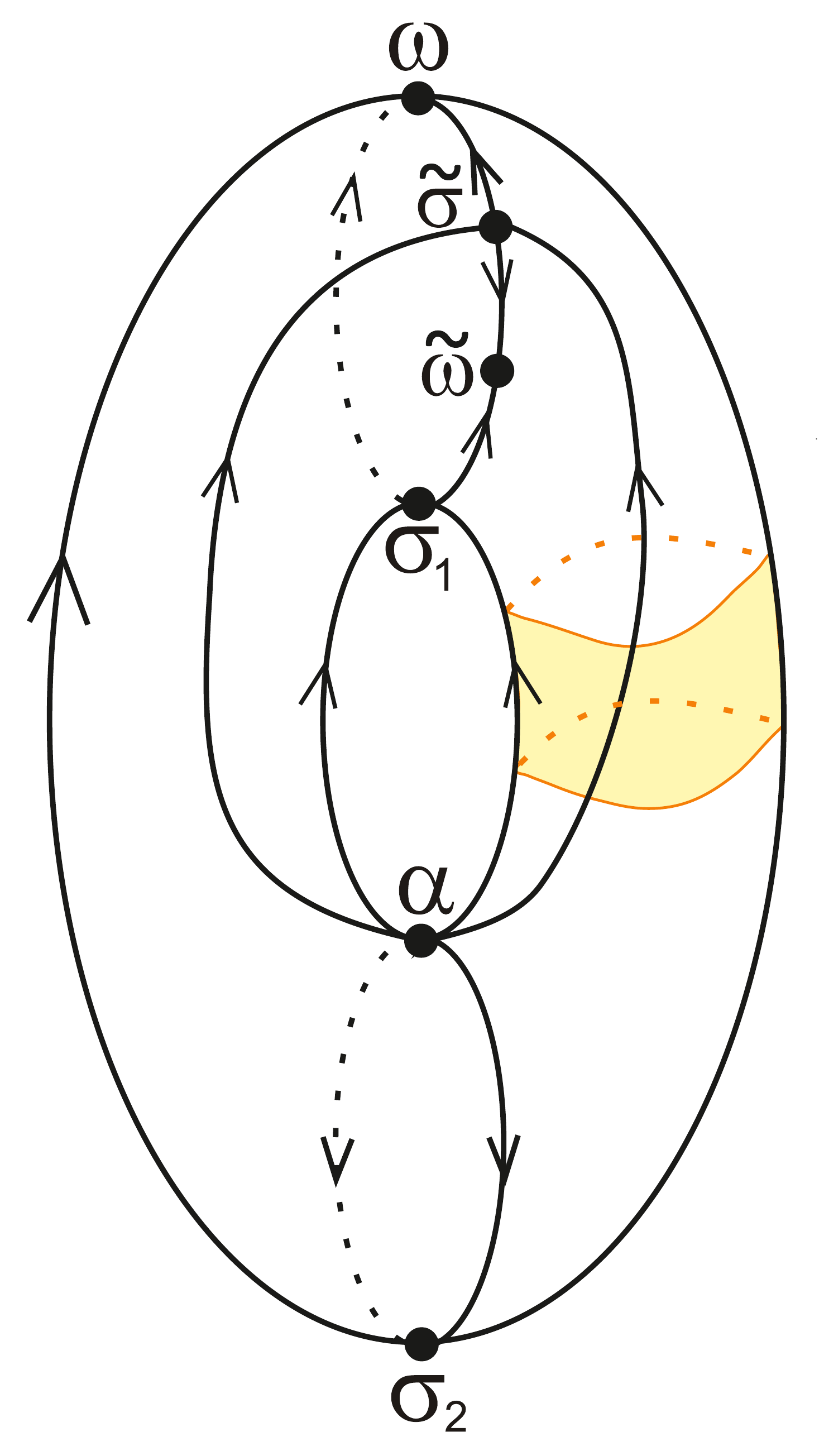} \\ а)$t=0$}
\end{minipage}
\hfill
\begin{minipage}[h]{0.32\linewidth}
\center{\includegraphics[width=0.9\linewidth]{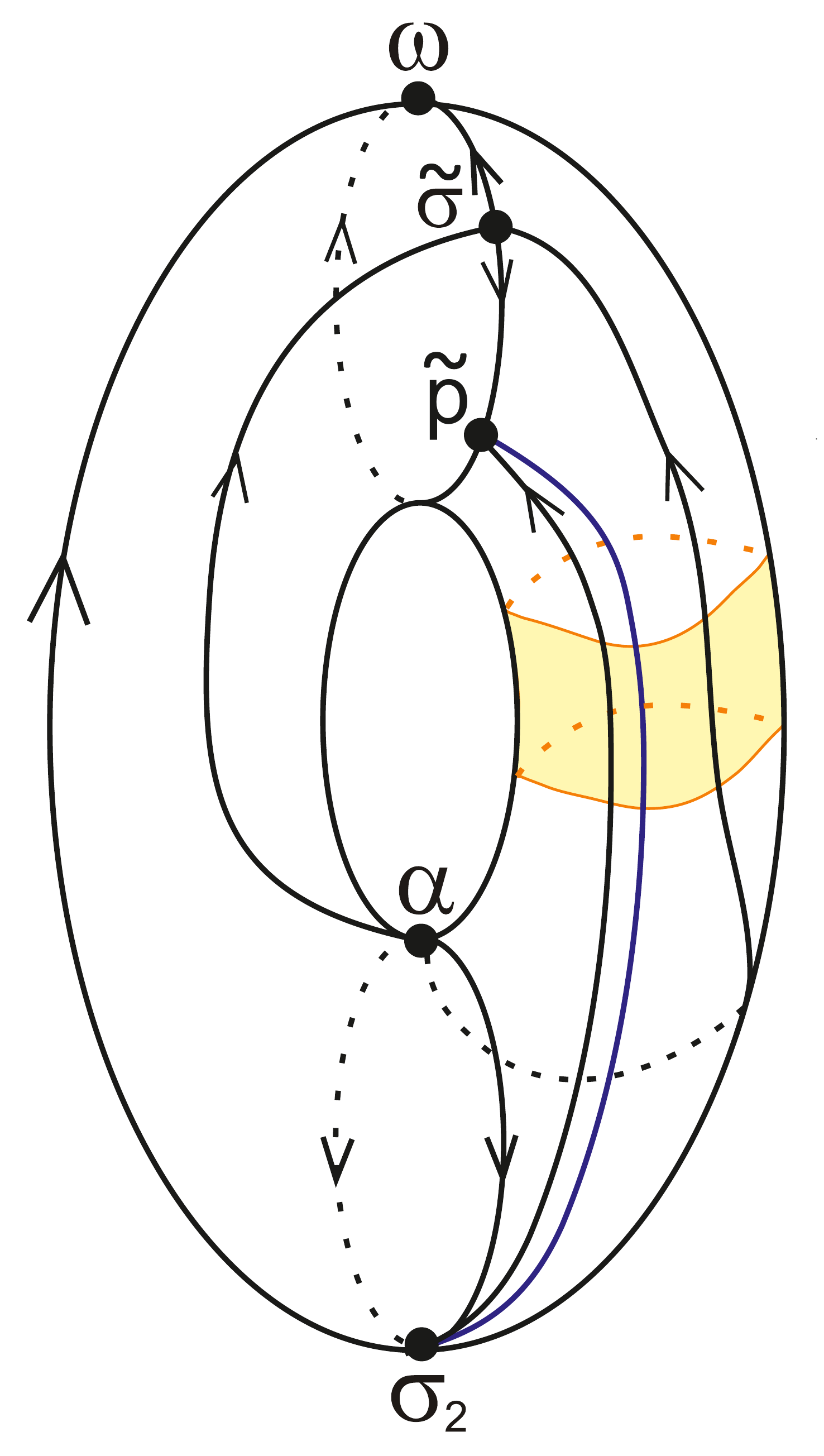} \\ б)$t=\frac12$}
\end{minipage}
\hfill 
\begin{minipage}[h]{0.32\linewidth}
\center{\includegraphics[width=0.9\linewidth]{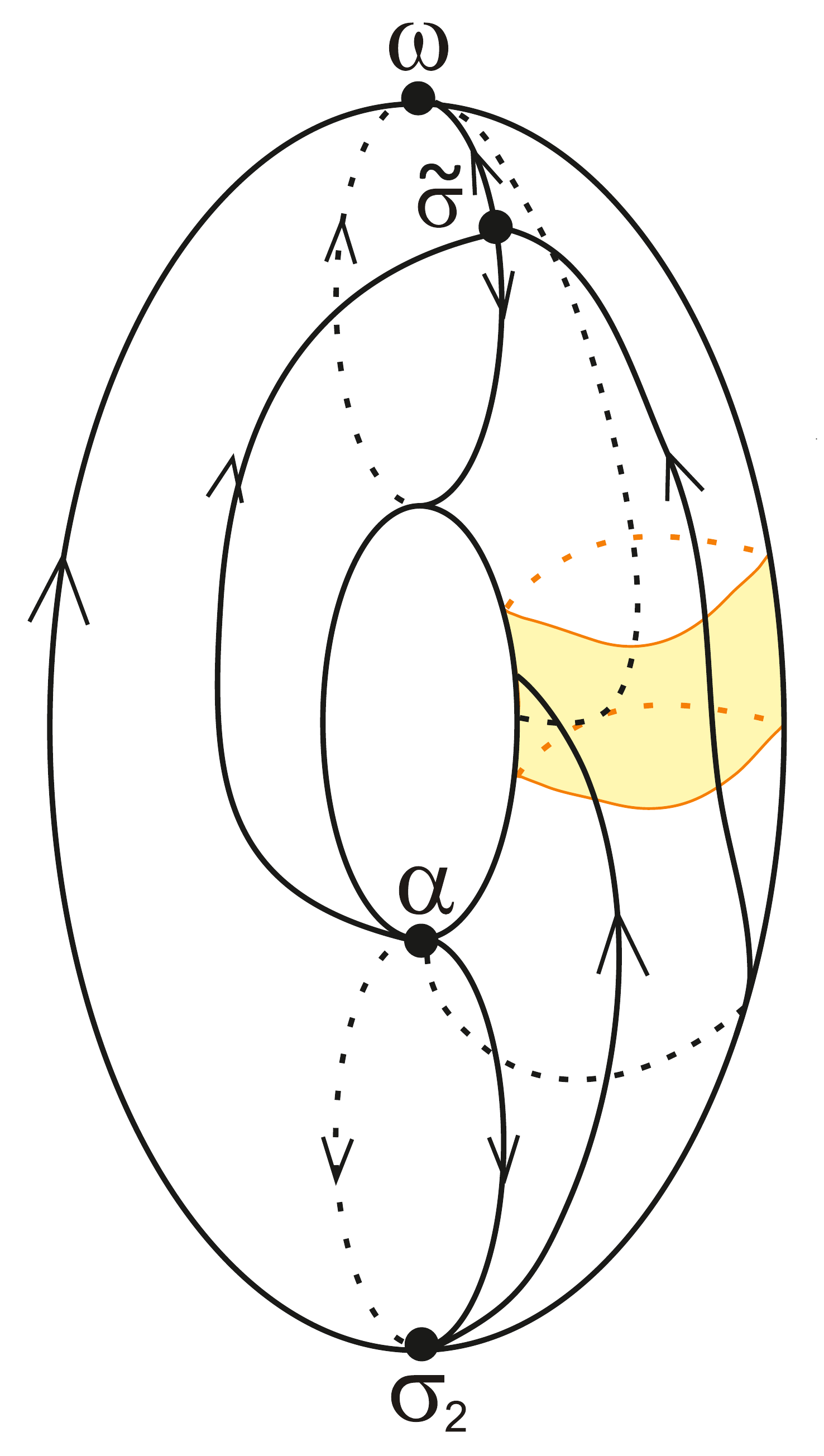} \\ б)$t=1$}
\end{minipage}
\hfill 
\caption{Isotopy $\Gamma^2_t$ on the torus.}
\label{ris:image1}
\end{figure}
According to lemma \ref{fA} the diffeomorphism $\Gamma^2_1$ can be connected by an arc without bifurcations $H_{\Gamma^2_1,t}$ with the diffeomorphism $f_{J_1}$.
\end{demo}

Denote by $H_{n,n+1,t}$ the arc with two saddle-node bifurcations connecting the diffeomorphisms $f_{J_n},\,f_{J_{n+1}}$ and given by the formula $$H_{n,n+1,t}=\hat J_n H_{0,1,t} \hat J_n^{-1}.$$ 

\subsection{Arc construction algorithm $H_{J,t}$}
In this section, using the model arcs constructed above, we will prove the following lemma.

\begin{lemma}\label{AA} 
The diffeomorphism $f_J$ is joined by a stable arc $H_{J,t}$ with a finite number of generically unfolding non-critical saddle-node bifurcations with the diffeomorphism $f_0$.

\end{lemma}
\begin{demo} Let $J= \begin{pmatrix}
  \mu^1& \mu^2\\
  \nu^1& \nu^2
\end{pmatrix}$ -- unimodular matrix such that $\mu^1\geq\mu^2\geq0$ and $\nu^1>\nu^2$, if $\mu^1=\mu^2$. Consider the following possibilities for the matrix $J$: 1) $\mu^2=0$; 2) $\mu^1=\mu^2=1$; 3) $\mu^2>\mu^1>0$. Construct the arc $H_{J,t}$ in each case separately.

In case 1) $J=J_n$. If $n>0$, then $H_{J_n,t}=H_{n-1,n,1-t}*\dots*H_{0,1,1-t}$ is the required arc. If $n<0$, then
$H_{J_n,t}=\hat{J}_{n} H_{J_{-n},1-t} \hat{J}_{n}^{-1}$ is the required arc. 

In case 2) $H_{J,t}=\hat{J} H_{J_{-1},1-t} \hat{J}^{-1}*H_{J_{\nu^2,t}}$ is the required arc.

In case 3) applying Euclid's algorithm to the pair $\mu_1,\mu_2$ generates a sequence of natural numbers $n_1,\dots,n_m,\,k_1,\dots,k_m$ тsuch that  $\mu^1=n_1 \mu^2 + k_1$, $\mu^2=n_2 k_1 + k_2$, $k_1=n_3 k_2 + k_3,\dots,k_{m-2}=n_m k_{m-1} + k_m$, где $k_{m-1}=1,\,k_m=0$. Let $k_{-1}=\mu^1,\,k_0=\mu^2$. Then the sequence $k_{-1},k_0,k_1,\dots,k_m$ satisfies the recurrence relation $$k_{i+1}=n_{i+1}k_i-k_{i-1},\,i=0, ... ,m-1.$$ Let $l_{-1}=\nu^1,\,l_0=\nu^2$ and define the sequence $l_{-1},l_0,l_1,\dots,l_m$ by the recurrent relation $$l_{i+1}=n_{i+1}l_i-l_{i-1},\,i=0, ... ,m-1.$$ Let $L_i=\begin{pmatrix}
  k_{i-1}& k_i\\
  l_{i-1}& l_i
\end{pmatrix}, i=0, ... ,m$. Then the arc $F_{i,t}=\hat L_{i-1}H_{J_{-n_i},t}\hat L_{i-1}^{-1},\,i=1, ... ,m$ joins diffeomorphisms $f_{L_{i-1}}$ and $f_{L_i}$ and contains $2n_i$ non-critical saddle-node bifurcations. Since $f_{L_m}=f_{J_{l_{m-1}}}$, then $H_{J,t}=F_{1,t}* ... *F_{m,t}*H_{J_{l_{m-1}}}$ is the required arc.
\end{demo}

\end{document}